\documentclass{amsart}
\usepackage[utf8]{inputenc}
\usepackage[T1]{fontenc}
\usepackage[english]{babel}
\usepackage[backend=biber, style=numeric, sorting=nty, maxbibnames=99, doi=false, isbn=false, eprint=false, giveninits=true]{biblatex}
\usepackage{csquotes}
\usepackage{latexsym}
\usepackage{amssymb}
\usepackage{dsfont}
\usepackage{mathtools}
\usepackage{mathrsfs}
\usepackage{enumitem}
\setlist{nolistsep}
\usepackage[hidelinks]{hyperref}
\usepackage{bm}
\usepackage{xcolor}

\allowdisplaybreaks[4]

\theoremstyle{plain}
\newtheorem{Theorem}{Theorem}[section]
\newtheorem*{Proposition*}{Proposition}
\newtheorem{Proposition}[Theorem]{Proposition}
\newtheorem{Lemma}[Theorem]{Lemma}
\newtheorem{Corollary}[Theorem]{Corollary}
\newtheorem{Definition}[Theorem]{Definition} 
\theoremstyle{definition}
\newtheorem{Remark}[Theorem]{Remark}	   

\numberwithin{equation}{section} 

\renewcommand{\div}{\operatorname{div}}
\DeclareMathOperator{\curl}{curl}

\newcommand{\R}{\mathbb{R}}

\newcommand{\N}{\mathbb{N}}

\newcommand{\abs}[1]{ \left\lvert#1\right\rvert} % absolute value: single vertical bars
\newcommand{\norm}[1]{\left\lVert#1\right\rVert} % norm: double vertical bars
\newcommand{\vertiii}[1]{{\left\vert\kern-0.25ex\left\vert\kern-0.25ex\left\vert #1 
		\right\vert\kern-0.25ex\right\vert\kern-0.25ex\right\vert}}
\newcommand{\supp}{\operatorname{supp}}

\let\oldexists\exists \let\exists\relax \DeclareMathOperator{\exists}{\oldexists} % fixes exist spacing
\let\oldforall\forall \let\forall\relax \DeclareMathOperator{\forall}{\oldforall}

\newcommand{\tr}{\operatorname{tr}}
\newcommand*\diff{\mathop{}\!\mathrm{d}}
\newcommand*\dd{\mathop{}\!\mathrm{d}}

\newcommand{\e}{\mathrm{e}}

\DeclareFieldFormat[unpublished]{title}{\emph{#1}}
\AtEveryBibitem{
	\ifentrytype{article}{}{\clearfield{pages}}
	\clearfield{edition}
	\clearfield{series}
	\ifentrytype{unpublished}{}{\clearfield{note}}
	\ifentrytype{book}{\clearfield{volume}}{}
	\ifentrytype{unpublished}{}{\clearfield{url}}
}
\renewcommand{\epsilon}{\varepsilon}
\newcommand{\Rsym}{\R^{3\times3}_{\mathrm{sym}}}
\newcommand{\epsilond}{\varepsilon^{\mathrm{d}}}
\newcommand{\mud}{\mu^{\mathrm{d}}}
\newcommand{\lambdad}{\lambda^{\mathrm{d}}}
\newcommand{\Tmax}{T_{\mathrm{max}}}

\begin{document}
\title[Decay of the Nonlinear Maxwell System with Absorbing Boundary]{Normal Trace Inequalities and Decay of Solutions to the Nonlinear Maxwell System with Absorbing Boundary}

\author{Richard Nutt}
\author{Roland Schnaubelt}
\address{RN, RS, Karlsruhe Institute of Technology\\
Department of Mathematics\\
Englerstraße~2\\
76131 Karlsruhe\\
Germany}
\email{richard.nutt@kit.edu}
\email{schnaubelt@kit.edu}

%\thanks{$^2$Corresponding author}

\subjclass[2020]{35Q61, 35L50, 35B40, 35B65}

\keywords{Quasilinear Maxwell system, absorbing boundary conditions, nonhomogenous anisotropic materials, trace regularity, exponential decay}

\begin{abstract}
	We study the quasilinear Maxwell system with a strictly positive, state dependent boundary conductivity.
	For small data we show that the solution exists for all times and decays exponentially to $0$.
	As in related literature we assume a nontrapping condition.
	Our approach relies on a new trace estimate for the corresponding non-autonomous linear problem, an observability-type estimate, and a detailed regularity analysis.
	The results are improved in the linear autonomous case, using properties of the Helmholtz decomposition in Sobolev spaces of (small) negative order.
\end{abstract}

\maketitle

\section{Introduction}

The Maxwell system is the foundation of electromagnetic theory.
It contains constitutive relations that describe the polarization $P$ and magnetization $M$ of the material in dependence of the electromagnetic fields.
In many physical models nonlinear effects occur which lead to nonlinear material laws, see e.g.\ \cite{agrawal13}, \cite{boyd2008}, \cite{butcher1990}, \cite{fabrizio2003}.
In this work we study instantaneous laws, see \cite{agrawal13} or \cite{fabrizio2003}, for which the Maxwell equations can be written as a quasilinear hyperbolic system.
It is well known that such systems can exhibit blow up, see e.g.\ \cite{dancona2018} in the Maxwell case.

We focus on the effect of a strictly positive surface conductivity $\lambda$, which may also depend on the electric field $E$.
In the case of linear material laws for $P$ and $M$, such nonlinear $\lambda$ had been studied in \cite{eller2002-1}, \cite{eller2002-2}, \cite{nicaise2003}, and for delayed problems in \cite{anikushyn2019}.
For state independent $\lambda$ and small initial fields it was shown in the paper \cite{pokojovy2020}, co-authored by one of us, that the solutions of the quasilinear Maxwell system~\eqref{eq:maxwell}--\eqref{eq:maxwell4} exist globally in time and decay exponentially as $t \to \infty$.
However, in \cite{pokojovy2020} it was assumed that the spatial domain $\Omega$ is strictly starshaped.
In related problems it is known that this assumption can be removed provided one can show extra regularity of the normal trace of the solutions, see \cite{lasiecka1992-2} for the wave equation and \cite{eller2007} for the linear autonomous Maxwell system with the boundary conditions of a perfect conductor (which are different from absorbing ones studied here).
In the present paper we prove a trace estimate for the solutions of linear non-autonomous anisotropic Maxwell systems with absorbing boundary conditions and inhomogeneities, see Corollary~\ref{cor:inequality} and Proposition~\ref{prop:gronwall}.
In Theorem~\ref{thm:normal-trace-vs-tang-trace} we improve the trace inequality in the autonomous case.
These results use properties of the Helmholtz decomposition and div-curl estimates also in negative order Sobolev spaces, which are shown in the appendix.

We use this extra trace regularity to establish an observability-type estimate for such systems in Proposition~\ref{prop:observability}. 
This result had been shown in \cite{pokojovy2020} for strictly starshaped domains only.
Based on this estimate, we can extend the analysis of \cite{pokojovy2020} to the case of a nonlinear boundary conductivity $\lambda(x,E)$ and show in our main Theorem~\ref{thm:mainThm} that solutions for small data converge exponentially to $0$ without assuming starshapedness.

We investigate the Maxwell system
\begin{align}
	\begin{split}\label{eq:maxwell}
		\partial_t (\epsilon(x,E(t,x)) E(t,x)) & = \curl H(t,x)\,, \\
		\partial_t (\mu(x,H(t,x)) H(t,x)) & = - \curl E(t,x)\,,
	\end{split}
	\qquad &&t\geq 0\,, x \in \Omega\,,\\
	\begin{split}\label{eq:maxwell2}
		\div(\epsilon(x,E(t,x))E(t,x)) &= 0\,,\\
		\div(\mu(x,H(t,x))H(t,x)) &= 0\,,
	\end{split}
	\qquad &&t\geq 0\,, x \in \Omega\,,\\
	H(t,x) \times \nu + \Big(\lambda\big(x,E(t,x)\times \nu\big)&(E(t,x) \times \nu) \Big) \times \nu = 0\,, &&t\geq 0\,, x \in \partial\Omega\,, \label{eq:maxwell3} \\
	E(0,x) = E^{(0)}(x)\,, &\quad H(0,x)= H^{(0)}(x)\,, && x \in \Omega\,,\label{eq:maxwell4}
\end{align}
on a bounded, smooth domain $\Omega \subseteq \R^3$ with connected complement, for the electric and magnetic fields $E(t,x), H(t,x) \in \R^3$ and given initial fields $(E^{(0)}, H^{(0)})$.
The permittivity $\epsilon(x,E)$, permeability $\mu(x,E)$, and surface conductivity $\lambda(x,E)$ may depend on position and state, and they belong to $\Rsym$.
So we have nonlinear, inhomogeneous, and anisotropic material laws.
As stated in Section~\ref{sec:2} we assume that the coefficients are $C^3$, symmetric, and uniformly positive definite at least for small $\abs{E}$, respectively $\abs{H}$.
In the analysis we often rewrite $\partial_t(\epsilon E)$ and $\partial_t(\mu H)$ as $\epsilond \partial_t E$ and $\mud \partial_t H$ for new coefficients $\epsilond$ and $\mud$ (see~\eqref{eq:epsilond}), which are supposed to have the same properties as $\epsilon$ and $\mu$.
This form of the equation facilitates energy estimates.

Local wellposedness of \eqref{eq:maxwell}--\eqref{eq:maxwell4} was shown in \cite{schnaubelt2021} for small data by energy methods.
(This smallness restriction is not necessary if $\lambda = \lambda(x)$ is state-independent.)
In this approach one has to control the Lipschitz norm of solutions, and thus their $H^3$-norms in the scale of integer-valued, $L^2$-based Sobolev spaces.
For this reason in \cite{schnaubelt2021} the nonlinear problem was solved in $H^m$ for $m \geq 3$.
(For full space problems one can reduce the necessary level of regularity below $\frac52$ by means of Strichartz estimates in some cases, cf.\ \cite{schippa2022}, but so far there are no such results for our boundary conditions.)
In this work we stick to $H^3$.
To bound the solutions in this norm, we look at the time-derived Maxwell systems \eqref{eq:inhomogenous_maxwell}--\eqref{eq:inhomogenous_maxwell_boundary}.
Here the coefficients $\epsilond$, $\mud$, and the analogue $\lambdad$ appear.
We stress that they are matrix-valued even if the given $\epsilon, \mu$, and $\lambda$ are scalar.
Moreover, they depend on time through the inserted solutions, and the system \eqref{eq:inhomogenous_maxwell}--\eqref{eq:inhomogenous_maxwell_boundary} contains error terms also at the boundary which will be treated as inhomogeneities.

To obtain $H^3$-solutions, the initial fields have to satisfy certain compatibility conditions stated in \eqref{eq:compatibility-condition}.
We note that these would simplify a lot for scalar-valued (isotropic) coefficients.
Applying the divergence to \eqref{eq:maxwell}, we see that the ``charges'' $\div(\epsilon E)$ and $\div(\mu H)$ are preserved in time.
We assume that the initial charges are $0$, see \eqref{eq:initial_solenoidality_conditions}, and that $\R^3 \setminus \Omega$ is connected in order to exclude non-zero stationary solutions of the form $(\nabla \varphi, \nabla \psi)$ where $\varphi$ and $\psi$ are constant on $\partial \Omega$, which would violate the desired decay property.

As in \cite{lasiecka2019} and \cite{schnaubelt2021}, the decay result follows from three propositions dealing with the time derivatives $\partial_t^k (E,H)$ for $k \in \{0,1,2,3\}$, see Section~\ref{sec:3}.
The fields $\partial_t^k (E,H)$ satisfy the boundary condition~\eqref{eq:maxwell3} up to lower-order terms.
This fact is crucial for the analysis.
An energy estimate and an observability-type estimate will allow us to control the squared $L^2$-norm of $\partial_t^k (E,H)$ by a dissipation term plus an error term which is small for small data, but contains space and time derivatives of higher order.
Surprisingly the regularity result Proposition~\ref{prop:regularity-boost} for the nonlinear problem allows us to absorb the error terms.
Theorem~\ref{thm:mainThm} on decay then follows by a standard bootstrap procedure, given at the end of Section~\ref{sec:3}.

The core observability-type estimate in Proposition~\ref{prop:observability} is based on a Morawetz multiplier argument as in \cite{pokojovy2020}, which uses ideas from \cite{eller2007} or \cite{nicaise2005} treating different boundary conditions.
For this result we have to assume the lower bound~\eqref{eq:technicalCond} on the radial derivatives of $\epsilon(x,0)$ and $\mu(x,0)$.
Heuristically this condition prevents trapping of the solution by back reflections so that they really reach the boundary where damping occurs.
The main difficulty is the control of boundary terms.
Tangential traces of solutions are bounded by the energy estimate, see also \cite{cagnol2011} or \cite{schnaubelt2021}.
To control also the normal trace, in Section~\ref{sec:4} we use the so-called collar operator studied in \cite{eller2007} in the context of different boundary conditions, cf.\ \cite{lasiecka1992-2} for earlier work on the wave equation.
This pseudodifferential operator allows us to trade space into time regularity in the course of a sophisticated regularity argument.
We further employ the div-curl estimate from Theorem~\ref{thm:div-curl-est} and exploit heavily the structure of the (time-differentiated) Maxwell system and the absorbing boundary condition. For Theorem~\ref{thm:div-curl-est} we have to assume that $\R^3\setminus\Omega$ is connected. 
In the trace estimates it is crucial that constants do not depend on time, see Corollary~\ref{cor:inequality}.
We also derive a more concise variant of the estimate in Proposition~\ref{prop:gronwall}, which has time dependent constants though.
This result was already shown  in \cite{cagnol2011} using completely different methods.
In Sections~\ref{sec:5} and \ref{sec:7} we then prove the observability-type estimate Proposition~\ref{prop:observability} and the regularity result Proposition~\ref{prop:regularity-boost}.

In Theorem~\ref{thm:normal-trace-vs-tang-trace} the trace estimate of Proposition~\ref{prop:gronwall} is improved for the linear, autonomous, anisotropic case~\eqref{pde:lin}.
We obtain a bound on the normal traces of $E$ and $H$ in $L^2((0,T) \times \partial \Omega)$ through the tangential trace of $E$, with constants being independent of the end time $T$.
As a by-product we also derive exponential stability of the system~\eqref{pde:lin}, which seems to be a new result for matrix-valued coefficients (and is not a special case of our nonlinear decay result).

In the proof of Theorem~\ref{thm:normal-trace-vs-tang-trace} we follow the approach of \cite{lasiecka1992-1} for the wave equation which involves a compactness argument to get rid of lower-order terms, cf.\ \cite{eller2007} or \cite{lasiecka1992-1}.
In order to perform this argument, we need the div-curl estimate from Theorem~\ref{thm:div-curl-est} in $H^{-\theta}(\Omega)$ for $\theta \in (0,\frac12)$, which requires several properties of the Helmholtz decomposition in $H^{-\theta}$ also established in the appendix.
Moreover, our reasoning involves two uniqueness results:
One for the evolution equation based on our observability estimates (see Lemma~\ref{lem:observability} and \eqref{eq:proof-time-independent-bnd-obs})
and a stationary one following from elliptic theory and the Helmholtz decomposition (see Theorem~\ref{thm:div-curl-est}~b)).

\section{Notation, assumptions and auxiliary results}
\label{sec:2}

As in \cite{pokojovy2020} we consider a bounded domain $\Omega \subseteq \R^3$ with $C^5$-boundary $\partial \Omega \eqqcolon \Gamma$.
% We introduce the time-space cylinders $\Omega_t \coloneqq (0,t)\times \Omega$, $\Gamma_t \coloneqq (0,t) \times \Gamma$, $\Omega_\infty \coloneqq \R\times \Omega$ and $\Gamma_\infty \coloneqq \R \times \Gamma$.
The outer unit normal will be denoted as $\nu \colon \Gamma \to \R^3$.
Furthermore, we introduce the function spaces
\begin{align*}
    C_\tau^k(\Gamma) &\coloneqq \{ f \in C^k(\Gamma) \mid f \cdot \nu = 0 \}\,,\\
    C^3_\tau(\Gamma \times \R^3, \Rsym) &\coloneqq \{A \in C^3(\Gamma \times \R^3, \Rsym) \mid A \nu^\perp \subseteq \nu^\perp \}\,,\\
	H^\theta(\div_\alpha) &\coloneqq \left\{f \in \big(H^\theta(\Omega)\big)^3 \mid \div (\alpha f) \eqqcolon \div_\alpha f \in H^\theta(\Omega)\right\}\,, \\
    H^\theta(\curl_\beta) &\coloneqq \left\{f \in \big(H^\theta(\Omega)\big)^3 \mid \curl(\beta f) \eqqcolon \curl_\beta f \in \big(H^\theta(\Omega)\big)^3\right\}\,,\\
    \intertext{for $\theta \in \R$ and suitable matrix-valued $\alpha,\beta \colon \Omega \to \R^{3 \times 3}$, and our solution space}
    G^k(\Omega) &\coloneqq G^k\left([0,\Tmax)\right) \coloneqq \bigcap_{j=1}^k C^k\big([0,\Tmax), H^{k-j}(\Omega)^6\big)\,.
\end{align*}
We also use the analogue of this space for compact time intervals, which is equipped with the canonical norm. For $\theta = 0$ we simply write $H(\div_\alpha) = H^0(\div_\alpha)$ and $H(\curl_\beta) = H^0(\curl_\beta)$, and if $\alpha=\beta = I_3$ is the identity matrix, we write $H^\theta(\div)$ and $H^\theta(\curl)$. Often we will drop $\Omega$ as well (e.g.\ set $G^k \coloneqq G^k(\Omega)$), and omit the power in $\big(L^2(\Omega)\big)^3$ etc. to shorten notation.
We assume that the permittivity $\epsilon$, the permeability $\mu$, and the surface conductivity $\lambda$ are of class
\begin{equation}\label{assumption:cont1}
  \epsilon, \mu \in C^3(\overline{\Omega}\times\R^3, \Rsym), \quad \lambda \in C^3_\tau(\Gamma \times \R^3, \Rsym)
\end{equation}
and are uniformly positive definite for zero fields, i.e., there is a constant $\eta$ such that
\begin{equation}\label{assumption:posdef1}
    \begin{aligned}
    \epsilon(x,0) \geq 2\eta I_3, \quad \mu(x,0) &\geq 2\eta I_3 \quad \text{for all $x \in \overline{\Omega}$} \\
    \text{and} \quad \lambda(x,0)&\geq 2 \eta I_3 \quad \text{for all $x\in \Gamma$}\,.
    \end{aligned}
\end{equation}
Note that the nonlinearities $\epsilon$ and $\mu$ are not required to be bounded in $E,H \in \R^3$.
Later on, however, the arguments $E$ and $H$ will be bounded and therefore so will be $\epsilon(\cdot,E)$ as well as $\mu(\cdot,H)$.

In a moment we will also discuss the time derived Maxwell equations for suitably smooth solutions.
In order to consider these systems, we introduce the matrices
\vspace{-5pt}
\begin{equation}\label{eq:epsilond}
    \begin{aligned}
    \epsilond_{ij}(x,\xi) &= \epsilon_{ij}(x,\xi) + \sum_{l=1}^3\partial_{\xi_j}\epsilon_{il}(x,\xi)\xi_l\qquad &&\text{for } x \in \overline{\Omega}, \xi \in \R^3,\\
    \mud_{ij}(x,\xi) &= \mu_{ij}(x,\xi) + \sum_{l=1}^3\partial_{\xi_j}\mu_{il}(x,\xi)\xi_l\qquad&&\text{for } x \in \overline{\Omega}, \xi \in \R^3,\\
    \lambdad_{ij}(x,\xi) &= \lambda_{ij}(x,\xi) + \sum_{l=1}^3\partial_{\xi_j}\lambda_{il}(x,\xi)\xi_l\qquad&&\text{for } x \in \Gamma, \xi \in \R^3,
    \end{aligned}
\end{equation}
where $i,j \in \{1,2,3\}$.
To unify notation, we set
\begin{align*}
    \widehat{\epsilon}_k &= \begin{cases}
        \epsilon(\cdot,E), \ &k=0,\\
        \epsilond(\cdot,E), \ &k\in\{1,2,3\},
    \end{cases}
    \,\qquad
    \widehat{\mu}_k = \begin{cases}
        \mu(\cdot,E), \ &k=0,\\
        \mud(\cdot,E), \ &k\in\{1,2,3\},
    \end{cases}
    \,\\
    \widehat{\lambda }_k &= \begin{cases}
        \lambda(\cdot,E), \ &k=0,\\
        \lambdad(\cdot,E), \ &k\in\{1,2,3\}.
    \end{cases}
\end{align*}
We also assume that
\begin{equation}
\begin{split}\label{assumption:cont2}
  &\partial_{\xi_j}\epsilon,   \partial_{\xi_j}\mu \in C^3(\overline{\Omega} \times \R^3, \R^{3 \times 3}), \quad \partial_{\xi_j}\lambda \in C^3(\Gamma\times \R^3, \R^{3 \times 3}), \\
  &\epsilond = (\epsilond)^\intercal,\quad \mud = (\mud)^\intercal,\quad \lambdad = (\lambdad)^\intercal\,.
\end{split}
\end{equation}
for $i,j \in \{1,2,3\}$.
Note that $\lambdad = \lambda$ if $\lambda = \lambda(x)$ is state-independent.
We extend $\lambda$ to a function on $\overline{\Omega} \times \R^3$ satisfying the same conditions as $\epsilon$ and $\mu$.
Since the coefficients are continuous, we also have uniform positivity at least for small fields; i.e.,
\begin{equation}\label{eq:uniformly-positive}
\begin{aligned}
	\epsilon(x,\xi), \epsilond(x,\xi) \geq \eta I\,, \quad \mu(x,\xi), \mud(x,\xi) \geq \eta I \quad \text{for all $\abs\xi \leq \delta_0, x \in \overline{\Omega}$}\,,\\
    \lambda(x,\xi), \lambdad(x,\xi) \geq \eta I \quad\text{for all $\abs\xi \leq \delta_0, x \in \Gamma$}\,,
\end{aligned}
\end{equation}
and some $\delta_0>0$.
An important special case is the Kerr law $\mu=\mu(x)$ and $\epsilon = \epsilon(x,\xi) = \epsilon_{\mathrm{lin}}(x) + \epsilon_{\mathrm{nl}}(x) \abs{\xi}^2$ for scalar coefficients with $\epsilon_{\mathrm{lin}} \geq 2 \eta$, cf.\ \cite{agrawal13} or \cite{fabrizio2003}.
Here one has $\epsilond(x, \xi) = \epsilon(x, \xi) + 2(\epsilon_{\mathrm{nl}}(x) \xi) \xi^\intercal$.
Anisotropic examples of polynomial type are discussed in Example~2.1 of \cite{lasiecka2019}.

Since the solutions of \eqref{eq:maxwell}--\eqref{eq:maxwell4} are supposed to satisfy the boundary condition at all times, it has to hold also for the initial values.
This leads to so-called ``compatibility conditions'' (of order 3) on $E^{(0)}$ and $H^{(0)}$.
Namely, for $E^{(0)}, H^{(0)} \in H^3(\Omega)$, a solution $(E,H) \in G^3$ of \eqref{eq:maxwell} possesses the time derivatives
\begin{equation}\label{eq:timeDerivativeAt0}
    \begin{split}
    E^{(1)} &\coloneqq \bigl(\epsilond(E^{(0)})\bigr)^{-1} \curl H^{(0)}\,,\\
    H^{(1)} &\coloneqq  - \bigl(\mud(H^{(0)})\bigr)^{-1} \curl E^{(0)}\,,\\
    E^{(2)} &\coloneqq  \bigl(\epsilond(E^{(0)})\bigr)^{-1} \left[\curl H^{(1)} - \Big(\sum_{\ell=1}^3 \partial_{\xi_\ell}\epsilond_{i,j}(E^{(0)})E^{(1)}_\ell \Big)_{i,j} E^{(1)} \right]\,,\\
    H^{(2)} &\coloneqq  - \bigl(\mud(H^{(0)})\bigr)^{-1}  \left[\curl E^{(1)} - \Big(\sum_{\ell=1}^3 \partial_{\xi_\ell}\mud_{i,j}(H^{(0)})H^{(1)}_\ell \Big)_{i,j} H^{(1)} \right]\,,
    \end{split}
\end{equation}
at time $0$, which leads to the compatibility conditions
\begin{equation}\label{eq:compatibility-condition}
\begin{aligned}
    &H^{(0)} \times \nu + \lambda(E^{(0)} \times \nu)(E^{(0)} \times \nu) \times \nu = 0\,,\\
    &H^{(1)} \times \nu + \lambdad(E^{(0)}\times \nu) (E^{(1)} \times \nu) \times \nu = 0\,, \\
    &\begin{aligned}
    H^{(2)} \times \nu &+ \lambdad(E^{(0)}\times \nu)(E^{(2)} \times \nu) \times \nu\\ &= -\bigg( \Bigl(\sum_{\ell=1}^3 \partial_{\xi_\ell} \lambdad_{i,j}(E^{(0)}\times \nu) (E^{(1)}\times \nu)_{\ell}\Bigr)_{i,j}E^{(1)}\times \nu\bigg) \times \nu\,,\end{aligned}
\end{aligned}
\end{equation}
on $\Gamma$.
Let $c_S$ be the norm of the Sobolev embedding $H^2(\Omega) \hookrightarrow C(\overline{\Omega})$.
Set $\tilde{\delta} \coloneqq \min\{1 , \frac{\delta_0}{c_S}\}$, cf.\ \eqref{eq:uniformly-positive}.
Under the above assumptions, in Theorem 6.4 of \cite{schnaubelt2021} it was shown that for $\delta \in (0, \tilde{\delta}]$ and sufficiently small initial values satisfying \eqref{eq:compatibility-condition}, i.e.,
\begin{equation}\label{eq:initial_data_bounded}
    \norm{(E^{(0)},H^{(0)})} _{H^3(\Omega)^6} \leq r (\delta) 
\end{equation}
for a constant $r(\delta)>0$ depending on $\delta$, there exists a unique classical solution $(E,H) \in G^3([0,\Tmax))$ such that $\Tmax > 1$ and 
\begin{equation}\label{eq:bound}
    \max_{0\leq j \leq 3} \left(\norm{\partial_t^jE(t)}_{H^{3-j}(\Omega)}^2 +\norm{\partial_t^j H(t)}_{H^{3-j}(\Omega)}^2 \right) \leq \delta^2 \quad \text{for } 0\leq t\leq T\,.
\end{equation}
We set 
\begin{equation}\label{eq:defTStar}
  T_* \coloneqq \sup \{\hat T \in [0,\Tmax] \mid \eqref{eq:bound} \text{ holds for } t \in [0,\hat T]\} \geq 1  \,.
\end{equation}
From the blow-up condition in Theorem 6.4 of \cite{schnaubelt2021} we infer that $T_* < \infty$ implies
\begin{equation}\label{eq:normBounded}
    \max_{0\leq j \leq 3}\left(\norm{\partial_t^jE(T_*)}_{H^{3-j}(\Omega)}^2 +\norm{\partial_t^j H(T_*)}_{H^{3-j}(\Omega)}^2 \right) = \delta^2 \,.
\end{equation}
Below we will assume $0 \leq t < T_*$ if we work with the solution $(E,H)$ to \eqref{eq:maxwell} -- \eqref{eq:maxwell4} with data satisfying \eqref{eq:initial_data_bounded}.

Furthermore, we define the commutator terms
\begin{align*}
     & \begin{aligned}
           f_0 & =f_1=0\,, & f_2 & = (\partial_t \epsilond(\cdot,E))\partial_t E\,,                                         & f_3 & =(\partial_t^2 \epsilond(\cdot,E))\partial_t E + 2 (\partial_t \epsilond(\cdot,E)) \partial_t^2 E\,, \\
           g_0 & =g_1=0\,, & g_2 & = (\partial_t \mud(\cdot,H))\partial_t H\,,                                                & g_3 & =(\partial_t^2 \mud(\cdot,H))\partial_t H + 2 (\partial_t \mud(\cdot,H)) \partial_t^2 H\,,           \\
           h_0 & =h_1=0\,,   & h_2 & = (\partial_t \lambdad(\cdot,E\times \nu))\partial_t E\times \nu\,, \hspace{-\linewidth} &     &
       \end{aligned} \\
     & h_3=(\partial_t^2 \lambdad(\cdot,E\times \nu))\partial_t E\times \nu + 2 (\partial_t \lambdad(\cdot,E\times \nu)) \partial_t^2 E\times \nu\,.
\end{align*}
Let $k \in \{0,1,2,3\}$.
With the definitions above the time-derived solutions satisfy
\begin{equation}\label{eq:derivative}
    \begin{aligned}
	\partial_t^k (\epsilon E ) &=  \widehat{\epsilon}_{k}\partial_t^{k} E + f_k\,, \qquad
	\partial_t^k (\mu H ) =  \widehat{\mu}_{k}\partial_t^{k} H + g_k \quad &&\text{on } \Omega\,,\\
    \partial_t^k (\lambda E\times\nu ) &= \widehat{\lambda}_k\partial_t^k E \times \nu + h_k \quad &&\text{on } \Gamma\,.
    \end{aligned}
\end{equation}
We thus arrive at the inhomogeneous system
\begin{align}
	&\begin{aligned}
        \partial_t (\widehat{\epsilon}_k \partial_t^k E) & =  \curl \partial_t^k H - \partial_t f_k\,,\\
		\partial_t (\widehat{\mu}_k \partial_t^k H) & = - \curl \partial_t^k E - \partial_t g_k\,,
	\end{aligned}
	\qquad t\geq 0, x \in \Omega, \label{eq:inhomogenous_maxwell}\\
    &\partial_t^k H \times \nu + \widehat{\lambda}_k(\partial_t^k E \times \nu) \times \nu = - h_k \times \nu \label{eq:inhomogenous_maxwell_boundary}\,.
\end{align}
Notice that \eqref{eq:maxwell}, \eqref{eq:maxwell2} and \eqref{eq:derivative} imply
\begin{equation}\label{eq:divergence_inhomogeneities}
  \div(\widehat{\epsilon}_k\partial_t^kE) = - \div(f_k)\, ,\quad   \div(\widehat{\mu}_k\partial_t^kH) = - \div(g_k)\,.
\end{equation}

Finally, in order to discuss the observability and energy estimates, we introduce the energy, dissipation, and error terms
\begin{align*}
	e_k(t) &= \tfrac12 \max_{0\leq j \leq k} \left(\norm{\widehat\epsilon_j^{1/2} \partial_t^jE(t)}_{L^2(\Omega)}^2 +\norm{\widehat\mu_j^{1/2} \partial_t^j H(t)}_{L^2(\Omega)}^2 \right), &e \coloneqq e_3\,,\\
	d_k(t) &= \max_{0\leq j \leq k} \norm{\lambda^{1/2} \tr_t \partial_t^jE(t)}_{L^2(\Gamma)}^2, &d \coloneqq d_3\,,\\
	z_k(t) &= \max_{0\leq j \leq k} \left(\norm{\partial_t^jE(t)}_{H^{k-j}(\Omega)}^2 +\norm{\partial_t^j H(t)}_{H^{k-j}(\Omega)}^2 \right), &z \coloneqq z_3\,,
\end{align*}
for the time-derived fields.
Here
\[
    \tr_t\colon H(\curl) \to \big(H^{-1/2}(\Gamma)\big)^3;\quad u \mapsto u \times \nu\,,
\]
denotes the tangential trace.
Similarly, we define the normal trace
\[
    \tr_n\colon H(\div) \to H^{-1/2}(\Gamma);\quad u \mapsto u \cdot \nu\,.
\]
These linear maps are bounded, see Theorem 2.2 and 2.3 of \cite{cessenat1996}.
We also use the rotated tangential trace
\[
    \tr_\tau \colon H(\curl) \to \big(H^{-1/2}(\Gamma)\big)^3;\quad u \mapsto \nu \times (u \times \nu) = \tr u - (\tr_n u)\nu\,.
\]
One can show the following estimates for the commutator terms
\begin{align}\label{eq:commutator-estimates}
    \begin{aligned}
	\max_{2 \leq k \leq 3, 0 \leq j \leq 1} \norm{\partial_t^j f_k}_{H^{4-j-k}(\Omega)} + \norm{\partial_t^j g_k}_{H^{4-j-k}(\Omega)} \lesssim z(t)\,,\\
    \max_{2 \leq k \leq 3} \norm{h_k}_{H^{3+1/2-k}(\Gamma)} \lesssim z(t)\,,
    \end{aligned}
\end{align}
see equation~(2.22) in \cite{pokojovy2020}.
We write ``$\lesssim$'' if the inequality holds up to a constant not depending on $t \in [0,T_*), T_*, \delta \in (0,\tilde{\delta}], r \in(0,r(\delta)]$, and $(E^{(0)}, H^{(0)})$ fulfilling \eqref{eq:compatibility-condition}, \eqref{eq:initial_data_bounded} and \eqref{eq:initial_solenoidality_conditions} below.
Moreover, such constants are denoted by $c,c_j,C$ or $C_j$.
We stress that $z$ is quadratic in the fields and bounded by $\delta z(t)^{1/2}$ due to \eqref{eq:bound}.

\begin{proof}[Sketch of proof of \eqref{eq:commutator-estimates}]
    We carry out the proof for $f_k$.
    For $g_k$ the argument can be repeated, and the same is true for $h_k$ after extending $h_k$ to $\overline{\Omega}$ using the trace theorem. 
    For $k \in \{0,1\}$ the commutator term $f_k$ vanishes and there is nothing to show.
    Notice that $H^2(\Omega) \hookrightarrow L^\infty(\Omega)$ and therefore $E$ and $\partial_t E$ are bounded.
    Since $\epsilond$ is assumed to belong to $C^3(\Omega \times \R^3)$, we can estimate the time derivative of 
    \begin{align*}
        f_3 &= \Big[\sum_{i,j=1}^3 \partial_{\xi_i,\xi_j} \epsilond(\cdot,E) \partial_t E_i \partial_t E_j\Big] \partial_t E+ \Big[\sum_{i=1}^3 \partial_{\xi_i}\epsilond(\cdot,E) \partial_t^2E_i \Big]\partial_t E \\
        &\hphantom{{}=} + 2\Big[\sum_{i=1}^3 \partial_{\xi_i} \epsilond(\cdot,E)\partial_t E_i\Big]\partial^2_t E 
        \eqqcolon f_{3,1} +f_{3,2} +2 f_{3,3}
    \end{align*}
    for instance by
    \begin{align*}
         \norm{\partial_t  f_{3,1}}_{L^2(\Omega)} &= \norm{\partial_t\Big[\Big(\sum_{i,j=1}^3 \partial_{\xi_i,\xi_j} \epsilond(\cdot,E) \partial_t E_i \partial_t E_j\Big) \partial_t E \Big]}_{L^2(\Omega)}\\
        &\leq \norm{\Big(\sum_{i,j, k=1}^3 \partial_{\xi_i,\xi_j, \xi_k} \epsilond(\cdot,E) \partial_t E_i \partial_t E_j\partial_t E_k\Big) \partial_t E}_{L^2(\Omega)} \\ 
        &\hphantom{{}=}+ 
        \norm{2 \Big(\!\sum_{i,j=1}^3 \partial_{\xi_i,\xi_j} \epsilond(\cdot,E) \partial^2_t E_i \partial_t E_j\Big) \partial_t E}_{L^2(\Omega)} 
        \\
        &\hphantom{{}=} + \norm{\Big(\!\sum_{i,j=1}^3 \partial_{\xi_i,\xi_j} \epsilond(\cdot,E) \partial_t E_i \partial_t E_j\Big) \partial^2_t E}_{L^2(\Omega)}\\
        &\lesssim \sum_{i,j=1}^3 \Big(\norm{\partial_t E_i}_{L^\infty(\Omega)}\norm{ \partial_t E_j}_{L^2(\Omega)} + \norm{\partial^2_t E_i}_{L^2(\Omega)}\norm{ \partial_t E_j}_{L^\infty(\Omega)} \\
        &\hphantom{{}=\sum_{i,j=1}^3 \Big(} + \norm{\partial_t E_i}_{L^\infty(\Omega)} \norm{\partial^2_t E_j}_{L^2(\Omega)} \Big)\\
        &\lesssim \sum_{i,j=1}^3 \Big( \norm{\partial_t E_i}_{H^2(\Omega)}\norm{ \partial_t E_j}_{L^2(\Omega)} + \norm{\partial^2_t E_i}_{L^2(\Omega)}\norm{ \partial_t E_j}_{H^2(\Omega)} \\
        &\hphantom{{}=\sum_{i,j=1}^3 \Big(} + \norm{\partial_t E_i}_{H^2(\Omega)} \norm{\partial^2_t E_j}_{L^2(\Omega)}\Big)\\
        &\lesssim z(t)\,.
    \end{align*}
    Analogous estimates for $\norm{\partial_t  f_{3,2}}_{L^2(\Omega)}$ and $\norm{\partial_t  f_{3,3}}_{L^2(\Omega)}$ can be shown, using also the embedding $H^1(\Omega) \hookrightarrow L^6(\Omega)$.
    In a similar way one obtains bounds on $ \norm{f_2}_{H^2(\Omega)} ,\norm{\partial_t f_2}_{H^1(\Omega)}$, and $\norm{f_3}_{H^1(\Omega)}$.
\end{proof}

We also define the function $m(x) = x-x_0$ on $\Omega$ for a fixed point $x_0 \in \Omega$.
It is used in the multiplier argument from \cite{pokojovy2020}.

\section{Main result for the nonlinear problem}
\label{sec:3}

We state our main decay theorem and the core ingredients for its proof.
Based on them, we show the theorem at the end of the section.
Note that compared to Theorem~2.2 in \cite{pokojovy2020} we no longer require the domain to be strictly starshaped and also permit semilinear boundary damping.
In \cite{pokojovy2020} the strict starshapedness of the domain was used to estimate trace terms, which we will treat here in a more delicate manner, utilizing methods from microlocal analysis as previously introduced in \cite{eller2007}.

We require the following \textit{non-trapping condition} on $\varepsilon$ and $\mu$, which emanates from the Morawetz multiplier technique used in \cite{pokojovy2020}:
\begin{equation}\label{eq:technicalCond}
	\begin{split}
		\epsilon(x,0)+\bigl(m(x) \cdot \nabla_x\bigr) \epsilon(x,0) &\geq \overline{\eta} \epsilon(x,0)\,,\\
		\mu(x,0)+\bigl(m(x) \cdot \nabla_x\bigr) \mu(x,0) &\geq \overline{\eta} \mu(x,0)
	\end{split}
\end{equation}
for a constant $\overline{\eta} >0$ and $x \in \overline{\Omega}$.
It says that $\epsilon$ and $\mu$ do not decay too rapidly in radial directions for small fields.
Heuristically this should reduce back reflections preventing the fields to reach the boundary.
Similar conditions were used in \cite{dancona2021}, \cite{eller2007}, and \cite{nicaise2005}, for instance.

\begin{Theorem}\label{thm:mainThm}
	Let $\Omega \subset \R^3$ be a bounded domain with $C^5$-boundary and connected complement.
	Furthermore, assume that the permittivity $\epsilon$ and permeability $\mu$ satisfy \eqref{assumption:cont1}, \eqref{assumption:posdef1}, \eqref{assumption:cont2}, and \eqref{eq:technicalCond}.
	We require that the initial values $E^{(0)}, H^{(0)} \in \left(H^3(\Omega)\right)^3$ satisfy the compatibility conditions \eqref{eq:compatibility-condition} as well as the initial ``charge'' conditions
	\begin{align}\label{eq:initial_solenoidality_conditions}
		\div\bigl(\epsilon(E^{(0)})E^{(0)}\bigr) = \div\bigl(\mu(H^{(0)})H^{(0)}\bigr) & = 0 \quad \text{on $\Omega$}\,.
	\end{align}
	Then there exist constants $M, \omega, r>0$ such that for $\norm{E^{(0)}}^2_{H^3(\Omega)} + \norm{H^{(0)}}^2_{H^3(\Omega)} \leq r^2$ the solutions $(E,H) \in G^3$ exist for all times $t \geq 0$, are unique, and decay exponentially, i.e.,
	\[
		\max_{0\leq j \leq 3} \norm{\left(\partial^j_t E(t), \partial^j_t H(t)\right)}^2_{H^{3-j}(\Omega)} 	\leq Me^{-\omega t} \norm{\left(E^{(0)}, H^{(0)}\right)}^2_{H^3(\Omega)}\,.
	\]
\end{Theorem}

Theorem \ref{thm:mainThm} is a consequence of the following three propositions on solutions to \eqref{eq:inhomogenous_maxwell} and \eqref{eq:inhomogenous_maxwell_boundary}: an energy inequality, an observability-type estimate and a regularity result.
All contain small error terms of highest order.
It is crucial that constants do not depend on time.

\begin{Proposition}\label{prop:energy_inequality}
	Assume the hypotheses of Theorem~\ref{thm:mainThm}, except for the connectedness of $\R^3\setminus\Omega$, \eqref{eq:technicalCond}, and \eqref{eq:initial_solenoidality_conditions}. Then we have
	\[
		e(t) + \int_s^t d(\tau) \dd \tau \leq e(s)+ c_1\int_s^t z^{3/2}(\tau) \dd\tau
	\]
	for $0 \leq s \leq t <T_*$.
\end{Proposition}

This proposition controls energy and dissipation by the initial energy plus error terms. The result was shown in Proposition~3.1 of \cite[]{pokojovy2020} for linear $\lambda = \lambda(x)$ and hence for $h_2=h_3=0$.
It is straightforward to extend the proof in \cite{pokojovy2020} to the present situation, using \eqref{eq:commutator-estimates} to estimate the arising terms with $h_2$ and $h_3$ by the $z^{3/2}$-integral.
In the next proposition we bound the time integral of the energy by the dissipation, the energy of initial and present time, as well as error terms.

\begin{Proposition}
	\label{prop:observability}
	Assume that the assumptions of Theorem~\ref{thm:mainThm} are satisfied. For $0 \leq s \leq t < T_*$ we obtain
	\[
		\int_{s}^{t} e(\tau) \diff \tau \leq c_2\int_s^t d(\tau) \diff \tau + c_3(e(t) + e(s)) +c_4 \int_s^t z^{3/2}(\tau) \diff \tau\,.
	\]
\end{Proposition}

This result is shown at the end of Section~\ref{sec:5}, and it improves Proposition~3.3 from \cite{pokojovy2020}.
We want to explain its role in the autonomous, homogenous and linear case.
In this setting the above estimates are true without $z$ and for all $t \geq 0$ and one has equality in Proposition~\ref{prop:energy_inequality} for $e_0, d_0$ and $z=0$ (see Lemma~3.2 in \cite{pokojovy2020}).
We thus obtain
\begin{equation}\label{eq:obs-est-linear}
\begin{aligned}
	T\Big(e_0(0) - \int_0^T d_0(\tau)\dd\tau\Big) = T e_0(T) &\leq \int_0^T e_0(\tau) \dd \tau \\
	 & \leq c_2 \int_0^T d_0(\tau) \dd \tau + c_3(e_0(T)+e_0(0)) \\
	 & \leq c_2 \int_0^T d_0(\tau) \dd \tau + 2c_3e_0(0)\,.
\end{aligned}
\end{equation}
Taking $T>2c_3$, we conclude
\begin{equation}\label{eq:initial_energy_bound_by_dissipation}
	(T-2c_3) e_0(0) \leq (c_2+T) \int_0^T d_0(\tau) \dd \tau\,.
\end{equation}
The dissipation thus controls the initial value, which is closely related to observation estimates.

Exactly as in Corollary~3.5 of \cite[]{pokojovy2020}, Propositions~\ref{prop:energy_inequality}~and~\ref{prop:observability} imply the estimate
\begin{equation}
	\label{eq:observability-like_inequality}
	e(t) + \int_s^t e(\tau) \dd \tau \leq C_1e(s) + C_2 \int_s^t z^{3/2}(\tau) \dd \tau
\end{equation}
for $0 \leq s \leq t < T_*$.
Without the $z$-term it is well known how to derive exponential stability from here. However, despite $\sqrt{z} \leq \delta$ being small, we cannot absorb the error term with $z$ by the left-hand side.
To achieve this, we have to bound the space derivatives of $(E,H)$ by their time derivatives.
This can be done exploiting the structure of the Maxwell systems \eqref{eq:maxwell}--\eqref{eq:maxwell4}, resulting in the following third ingredient of the proof of Theorem~\ref{thm:mainThm}.

%To reference the following Proposition later on
\edef\regBoostThmCnt{\arabic{Theorem}}
\edef\regBoostSecCnt{\arabic{section}}

\begin{Proposition}\label{prop:regularity-boost}
	Under the assumptions of Theorem~\ref{thm:mainThm}, with the exception of \eqref{eq:technicalCond} and the connectedness of $\R^3\setminus\Omega$, there exist constants $c_5, c_6>0$ such that
	\[
		z(t) \leq c_5e(t) + c_6z^2(t)	
	\]
	for $0 \leq t < T_*$.
\end{Proposition}

We sketch the proof of this result in Section~\ref{sec:7}, heavily relying on the proof of Proposition~4.1 in \cite{pokojovy2020}.
Based on the propositions above, we can now show our first main result.

\begin{proof}[Proof of Theorem~\ref{thm:mainThm}]
	The reasoning follows the lines of the proof of Theorem~2.2 of \cite[]{pokojovy2020}.
	For convenience, we present the arguments briefly. We first combine estimate \eqref{eq:observability-like_inequality} and Proposition~\ref{prop:regularity-boost}.
	Fixing a sufficiently small $\delta = \overline{\delta} \in (0, \tilde{\delta}]$, we derive
	\begin{equation}\label{eq:regularity-boost1}
		z(t) + \int_s^t z(\tau) \dd \tau \leq C z(s)
	\end{equation}
	for all $0 \leq s \leq t <T_*$ and a constant $C \geq 1$.
	The radius $\tilde\delta$ was defined before \eqref{eq:initial_data_bounded}.
	Suppose that $T_*<\infty$.
	It follows $z(T_*) = \overline{\delta}^2$ by \eqref{eq:normBounded}.
	On the other hand, inequality~\eqref{eq:regularity-boost1} yields $z(t) \leq Cz(s)$. For initial fields with $\norm{\left(E^{(0)}, H^{(0)}\right)}_{H^3} \leq r$, one can bound $z(0) \leq c_0r^2$ using equations \eqref{eq:timeDerivativeAt0}.
	As a result for a fixed, sufficiently small radius $\overline{r} \in (0,r(\overline{\delta})]$ we infer the contradiction $z(T_*) \leq \frac12 \overline{\delta}^2$.
	Therefore, the solution exists and \eqref{eq:regularity-boost1} is valid for all times.
	For $T = C(2C-1)$ we deduce
	\begin{align*}
		z(T)+ \frac{T}{C} z(T) &\leq z(T)+ \int_0^T z(\tau) \dd \tau \leq C z(0)\,,\\
		z(T) &\leq \frac{C}{1+ \frac{T}{C}} z(0) = \frac12 z(0)\,.
	\end{align*}
	The assertion follows by iteration.
\end{proof}

	Before proving Proposition~\ref{prop:observability} in the next two sections, we simplify the notation. Instead of \eqref{eq:inhomogenous_maxwell} and \eqref{eq:inhomogenous_maxwell_boundary} we study the linear non-autonomous system
	\begin{equation} \label{pde:1}
		\begin{aligned}
		\begin{aligned}
			\partial_t (\alpha u) & = \curl v + \partial_t \varphi\,, \\
			\partial_t (\beta v) & = - \curl u + \partial_t\psi\,, 
		\end{aligned} 
		\hspace{-79.5pt}&\hspace{79.5pt}\qquad &&t\geq 0\,, x \in \Omega\,, \\ %alignment hack, better way?
			v \times \nu + (\gamma(u\times\nu))\times \nu &= \omega\,, \qquad&& t\geq 0\,, x \in \Gamma\,,\\
			u(0) &= u^{(0)}\,, \quad v(0)=v^{(0)}\,, \qquad&& x \in \Omega\,.\\
		\end{aligned}
	\end{equation}
	Where we assume that
	\begin{align}\label{eq:regularityCoefficients}
		\begin{split}
			\alpha, \beta &\in C^1\bigl([0,T] \times \overline{\Omega}, \R^{3 \times 3}_\text{sym}\bigr), \quad \gamma \in C_\tau^1\bigl([0,T] \times \Gamma, \R^{3 \times 3}_\text{sym}\bigr), \quad\text{with }\alpha, \beta, \gamma \geq \eta\,,\\
			\varphi, \psi &\in G^1 = C^1\big([0,T], L^2(\Omega)^3\big) \cap C\big([0,T], H^1(\Omega)^3\big)\\
			\omega &\in C([0,T], H^{1/2}(\Gamma)) \text{ with } \nu \cdot \omega = 0\,,
		\end{split}
	\end{align}
	and choose initial values $u^{(0)}, v^{(0)} \in (L^2(\Omega))^3$ such that
	\begin{equation}\label{eq:regularityCoefficients2}
		 \div\big(\alpha(0)u^{(0)}\big) = \div \varphi(0) \text{ and } \div\big(\beta(0)v^{(0)}\big) = \div \psi(0)\,.
	\end{equation}
	To recover the original system \eqref{eq:inhomogenous_maxwell} and \eqref{eq:inhomogenous_maxwell_boundary}, we simply resubstitute $\alpha = \widehat{\epsilon}_k, \beta = \widehat{\mu}_k, \gamma = \widehat{\lambda}_k, \varphi = -f_k, \psi = -g_k$ and $\omega = -h_k\times \nu$ for $k \in \{0,1,2,3\}$.

According to Proposition~3.1 of \cite{schnaubelt2021} there exists a unique weak solution $(u,v) \in C^0([0,T], (L^2(\Omega))^6)$ of \eqref{pde:1} with tangential trace $(\operatorname{tr}_\tau u,\operatorname{tr}_\tau v) \in L^2(\Omega_T, \R^6)$. Observe that $\div(\alpha u) = \div(\varphi)$ and $\div(\beta v) = \div(\psi)$ belong to $C([0,T], L^2(\Omega))$.

\section{Trace inequality}
\label{sec:4}

The main goal of this section is to establish new bounds on the normal traces of the electric and magnetic fields and their derivatives, which are needed to show Proposition~\ref{prop:observability}.
We do not assume \eqref{eq:technicalCond} and that $\R^3\setminus\Omega$ is connected in this section.

\subsection*{The pseudodifferential collar operator \texorpdfstring{$X$}{X}}

To obtain the desired trace regularity, we introduce the pseudodifferential operator $X$ as defined in Definition 2.1 of \cite{eller2007}.
This operator has the advantageous property that it allows us to trade time regularity for spatial regularity on the boundary $\Gamma$ of our domain, as shown in Lemma~2.2 of \cite{eller2007}.
We restate its definition.

\begin{Definition} \label{def:pdo}
	Since our domain $\Omega$ has, in particular, a $C^2$ boundary there exists a tubular neighborhood $U$ of $\Omega$, and a partition $\Omega \cup U = \bigcup_{i=0}^m U_i$ by open $U_i\subseteq \R^3$ such that
	\begin{itemize}[label=--]
		\item $U_0 \subseteq \Omega$,
		\item $U_i \subseteq U$ with $U_i \cap \partial \Omega \neq \emptyset$ for $i \geq 1$, and
		\item there exists a partition of unity $(\phi_i)_{i=0,\dots,m}$ on $\Omega$ such that $\phi_i \in C_c^\infty(U_i)$ and $\sum_{i} \phi_i (x)= 1$ for $x \in \Omega$.
	\end{itemize}
	Furthermore, we choose a neighborhood $V_i \subseteq \R^3$ around the origin and coordinate mappings $\Phi_i\colon V_i \to U_i \in C^2(V_i)$ such that
	\begin{itemize}[label=--]
		\item $\Phi_i(V_i \cap \{(x,y,z) \mid z \geq 0\}) = U_i \cap \Omega$,
		\item $\Phi_i(V_i \cap \{(x,y,z) \mid z = 0\}) = U_i \cap \partial\Omega$.
	\end{itemize}
	We begin constructing $X$ locally on $\R \times U_i$.
	Let $\varkappa>0$, $\tilde{\chi} \in C^\infty(\R^3\setminus\{0\})$ with $0 \leq \tilde{\chi} \leq 1$ and
	\[
		\tilde{\chi}(\xi_0, \eta_1, \eta_2) = \begin{cases}
			1\,, \qquad & \abs{\xi_0} \leq \varkappa \abs{(\eta_1, \eta_2)}/2\,, \\
			0\,, \qquad & \abs{\xi_0} \geq \varkappa \abs{(\eta_1, \eta_2)}\,,
		\end{cases}
	\]
	and $0 \leq \rho \leq 1$ be a cutoff at the origin with
	\[
		\rho(\xi_0, \eta_1, \eta_2) = \begin{cases}
			1\,, \qquad & \abs{(\xi_0, \eta_1, \eta_2)}^2 \geq 2\,, \\
			0\,, \qquad & \abs{(\xi_0, \eta_1, \eta_2)}^2 \leq 1\,.
		\end{cases}
	\]
	We set $\chi(\xi_0, \eta_1, \eta_2, \eta_3) = \left(\rho\tilde{\chi}\right)(\xi_0, \eta_1, \eta_2) \in C^\infty(\R^4)$.
	For $u \in L^2(\R, H^s(\Omega))$ define $u_i = \phi_i u$ and $v_i(\xi_0,\eta) = u_i(\xi_0,\Phi_i(\eta))$ for $(\xi_0,\eta) \in \R \times V_i$.
	Take $\sigma_i \in C^\infty_0(V_i)$ such that $\sigma_i|_{\supp(\phi_i \circ \Phi_i)} = 1 $.
	The operator $X$ is then given by
	\[
		X(x,D)u = \sum_{i=1}^m [\sigma_i \chi(D_t, D_{y_1}, D_{y_2})v_i] \circ \Phi_i^{-1}\,.
	\]
	Its symbol is of class $C^1S_{cl}^0$, cf.\ \cite{taylor1991} and \cite{eller2007}.

	The restriction of $X$ to the boundary $\Gamma$ (as a pseudodifferential operator on a manifold) is also denoted by $X$.
\end{Definition}

We recall the central properties of this operator as stated and proven in Lemmas~2.2 and 2.3 of \cite{eller2007}, setting $\Gamma_\infty \coloneqq \R \times \Gamma$.
Observe that in the next Lemma the constant $K_1$ can be chosen arbitrarily small by fixing a small $\varkappa>0$ in Definition~\ref{def:pdo}.

\begin{Lemma}\label{lemma:eller:2.2}
	Let $u \in L^2\big(\R, H^s(\Omega)\big)$. Then there exists a constant $K_1>0$ such that
	\[
		\int_\R \norm{X\partial_t u}_{H^{s-1}(\Omega)}^2 \dd \tau \leq K_1 \left[ \int_\R \norm{X u}_{H^{s}(\Omega)}^2 \dd \tau + \int_\R \norm{u}_{H^{s-1}(\Omega)}^2 \dd \tau\right]
	\]
	for $0\leq s \leq 1$.
	The constant $K_1=K_1(\varkappa)$ belongs to $\mathcal{O}(\varkappa^2)$ for $\varkappa \to 0$.
\end{Lemma}

\begin{Lemma}\label{lemma:eller:2.3}
	Let $u \in L^2(\Gamma_\infty)$. Then there exists a constant $K_2=K_2(\varkappa)>0$ with
	\[
		\norm{(1-X) u}_{L^2(\Gamma_\infty)}^2 \leq K_2 \left[ \norm{u}_{H^{-1}(\Gamma_\infty)}^2 + \norm{\partial_t u}_{H^{-1}(\Gamma_\infty)}^2 \right]\,.
	\]
\end{Lemma}

Using the fact that $X$ is essentially a Fourier-multiplier at least locally and Chapter~4 of \cite{taylor1991}, we obtain the following regularity result.

\begin{Lemma}\label{lemma:X_cont}
	Let $-1\leq s\leq 1$. Then $X \colon L^2(\R,H^s(\Omega)) \to L^2(\R,H^s(\Omega))$ and $X \colon L^2(\R, H^s(\Gamma)) \to L^2(\R, H^s(\Gamma))$ are continuous.
\end{Lemma}

\subsection*{Trace inequality}

For a moment we consider the slightly more general inhomogeneities $\xi, \zeta \in L^2([0,T] \times \Omega)$ on the right-hand side of \eqref{pde:1}, i.e., the system
\begin{equation} \label{pde:2}
	\begin{aligned}
		\partial_t (\alpha(t,x) u(t,x)) & = \curl v(t,x) + \xi\,,\qquad &&t\geq 0\,, x \in \Omega\,, \\
		\partial_t (\beta(t,x) v(t,x)) & = - \curl u(t,x) + \zeta\,, \qquad&& t\geq 0\,, x \in \Omega\,, \\
		v \times \nu + (\gamma(u\times\nu))\nu & = \omega\,, \qquad &&t\geq 0\,, x \in \Gamma,\\
		u(0) &= u^{(0)}\,, \quad v(0)=v^{(0)}\,, \\
	\end{aligned}
\end{equation}
with coefficients and data as in \eqref{eq:regularityCoefficients} and \eqref{eq:regularityCoefficients2} above.
Then Proposition 3.1 of \cite{schnaubelt2021} still yields a unique weak solution $(u,v) \in C^0([0,T], (L^2(\Omega))^6)$.

In order to apply the pseudodifferential operator $X$ from Definition \ref{def:pdo} to the solutions $(u,v)$ defined on $[0,T]$, we have to extend $(u,v)$ by $0$ to $L^2(\R, L^2(\Omega))$.
We still denote this extension by $u$ respectively $v$, and proceed analogously with $\xi ,\zeta$ and $\omega$.
The following div-curl estimate allows us to bound the image of the pseudodifferential operator in $H^1$.
Here and below we often omit variables.
Recall that we also write $L^2$ instead of $\big(L^2(\Omega)\big)^3$, etc.
Initially we assume that $(u,v)$ belong to $G^1=G^1([0,T])$.
This assumption will be removed in Corollary~\ref{cor:inequality}.
We note that below the case $\theta = 0$ suffices for the proof of Proposition~\ref{prop:observability}. Negative $\theta$ are needed in Section~\ref{sec:6} for improved results on the autonomous linear problem.

\begin{Lemma}
	\label{lem:DivCurl}
	Under the above assumptions, let $(u,v) \in G^1$ be solutions to the system \eqref{pde:2} and $\theta \in [0,\frac12)$.
	We then have
	\begin{multline*}
		\int_\R \norm{X (u,v)}^2_{H^{1-\theta}(\Omega)} \dd \tau \lesssim  \int_\R \Big( \norm{X \div(\alpha u, \beta v)}^2_{H^{-\theta}(\Omega)} 
		+ \norm{X \curl (u,v)}^2_{H^{-\theta}(\Omega)} \\
		+ \norm{ (u,v)}^2_{H^{-\theta}(\Omega)} + \norm{ (u,v)}^2_{H^{-1/2-\theta}(\Gamma)} + \norm{ \omega}^2_{H^{1/2-\theta}(\Gamma)} \Big)\dd \tau.
	\end{multline*}
\end{Lemma}

\begin{proof}
	By Lemma~\ref{lemma:X_cont} the operator $X \colon L^2(\R,H^s(\Omega)) \to L^2(\R,H^s(\Omega))$ is continuous for $-1\leq s\leq 1$.
	We apply Theorem~\ref{thm:div-curl-est}~a) to $X u$ and $X v$, obtaining
	\begin{align} \label{eq:DivCurl}
				\notag	\norm{X (u,v)}_{H^{1-\theta}(\Omega)} &\lesssim \norm{ \div_\alpha X u}_{H^{-\theta}(\Omega)} + \norm{ \div_\beta X u}_{H^{-\theta}(\Omega)}  + \norm{X(u,v)}_{H^{-\theta}(\Omega)} \\
					&\hphantom{={}}+ \norm{ \curl X (u,v)}_{H^{-\theta}(\Omega)}\\ 
				\notag	&\hphantom{={}}+ \norm{X v \times \nu + (\gamma((X u)\times\nu))\times\nu}_{H^{1/2-\theta}(\Gamma)}\,.
	\end{align}
	% Note that there is a typo in second line of the proof in Lemma~5.1 of \cite{pokojovy2020}. It should be $\norm{w}_{H^1} \lesssim \norm{\curl u}_{L^2}$.
	To estimate the last term, we note that the cross product with $\nu$ can be written as a multiplication with a matrix $B$ whose entries and their derivatives are bounded.
	The same holds for $\gamma$ and thus for the entries of $B\gamma B$.
	We write
	\[
		(X v )\times \nu + (\gamma((X u)\times\nu))\times\nu = B (X v) + B\gamma B (X u)
	\]
	The commutators $[X,B],[X,B\gamma B] \colon L^2(\R, H^{-s}(\Gamma)) \to L^2(\R,H^{1-s}(\Gamma))$ are continuous for $s\in [0,1]$ (see Proposition~4.1.E of \cite{taylor1991}).
	We then compute
	\begin{align*}
		\int_\R & \norm{X v \times \nu + (\gamma((X u)\times\nu))\times\nu}^2_{H^{1/2-\theta}(\Gamma)} \dd \tau \\
		 & = \int_\R\norm{X( v \times \nu + \gamma( u\times\nu)\times\nu) + [X,B] v + [X,B\gamma B] u}_{H^{1/2-\theta}(\Gamma)}^2 \dd \tau \\
		 & \lesssim \int_\R\left(\norm{ \omega}^2_{H^{1/2-\theta}(\Gamma)}+ \norm{(u,v)}^2_{H^{-1/2-\theta}(\Gamma)} \right)\dd \tau\,,
	\end{align*}
	using the boundary condition $ v \times \nu + \gamma( u\times\nu)\times\nu = {w}$ of \eqref{pde:2} and the continuity of $X$ on $\Gamma$.
	Since $\curl$, $\div_\alpha$ and $\div_\beta$ have $C^1S^1_{cl}$-symbols, the commutators with $X$ are bounded on $H^{-\theta}$ which can be shown as in Proposition 4.1.E of \cite{taylor1991}.
	The assertion now easily follows.
\end{proof}

Exploiting the smallness of $K_1(\varkappa)$ in Lemma~\ref{lemma:eller:2.2}, we can simplify the estimate above.

\begin{Corollary}\label{cor:DivCurl2}
	In the situation of Lemma~\ref{lem:DivCurl} there is a number $\varkappa_0>0$ so that for all $\varkappa \in (0, \varkappa_0]$ in Definition~\ref{def:pdo} we have
	\begin{multline*}
		\int_\R \norm{X (u,v)}^2_{H^{1-\theta}(\Omega)}  \dd \tau
		\lesssim \int_\R  \Big(\norm{(u,v)}^2_{H^{-\theta}(\div_\alpha)\times H^{-\theta}(\div_\beta)}\\ 
		+ \norm{(u,v)}^2_{H^{-1/2-\theta}(\Gamma)} 
		+ \norm{X( \xi, \zeta)}^2_{H^{-\theta}(\Omega)} + \norm{ \omega}_{H^{1/2-\theta}(\Gamma)}^2 \Big)\dd \tau \,.
	\end{multline*}
\end{Corollary}

\begin{proof}
	Via \eqref{pde:2} we rewrite the curl-terms in Lemma~\ref{lem:DivCurl} as
	\[
		\norm{X \curl u}_{H^{-\theta}}^2 + \norm{X \curl v}_{H^{-\theta}}^2 = \norm{X \partial_t (\alpha u) - X \xi}_{H^{-\theta}}^2 + \norm{X \partial_t(\beta v) - X \zeta}_{H^{-\theta}}^2\,.
	\]
	By Chapter~4 in \cite[]{taylor1991} the commutator $[X,\alpha] \colon L^2(\R,H^{-\theta}(\Omega)) \to L^2(\R,H^{1-\theta}(\Omega))$ is bounded.
	This fact and Lemmas~\ref{lemma:eller:2.2}, \ref{lemma:X_cont}, and \ref{lem:DivCurl} lead to
	\begin{align*}
		\int_\R& \norm{X (u,v)}_{H^{1-\theta}}^2\dd \tau  \lesssim \int_\R \Big(\norm{\div(\alpha u,\beta v)}_{H^{-\theta}}^2 + K_1(\varkappa)\norm{X (\alpha u, \beta v)}_{H^{1-\theta}}^2\\
		& \quad + \norm{X(\xi,\zeta)}^2_{H^{-\theta}} + \norm{(u,v)}_{H^{-\theta}}^2 + \norm{(u,v)}_{H^{-1/2-\theta}(\Gamma)}^2  + \norm{\omega}^2_{H^{1/2-\theta}(\Gamma)} \Big)\dd \tau \\
		 & \lesssim \int_\R \Big(\norm{\div(\alpha u, \beta v)}^2_{H^{-\theta}} + K_1(\varkappa)\norm{X (u,v)}_{H^{1-\theta}}^2 + \norm{(u,v)}_{H^{-\theta}}^2\\
		 & \quad + \norm{(u,v)}_{H^{-1/2-\theta}(\Gamma)}^2 + \norm{X(\xi,\zeta)}^2_{H^{-\theta}}+ \norm{\omega}^2_{H^{1/2-\theta}(\Gamma)}\Big)\dd \tau\,,
	\end{align*}
	Choosing $\varkappa$ and hence $K_1(\varkappa)$ small enough yields
	\begin{align*}
		\int_\R\norm{X (u,v)}_{H^{1-\theta}}^2\dd \tau  &\lesssim \int_\R \Big(\norm{u}^2_{H^{-\theta}(\div_\alpha,\Omega)} + \norm{ v}^2_{H^{-\theta}(\div_\beta,\Omega)} + \norm{(u,v)}^2_{H^{-1/2-\theta}(\Gamma)}\\
		 & \quad +\norm{X( \xi, \zeta)}^2_{H^{-\theta}} + \norm{\omega}^2_{H^{1/2-\theta}(\Gamma)} \Big)\dd \tau \,. \qedhere
	\end{align*}
\end{proof}

From now on $\varkappa=\varkappa_0$ is fixed in the definition of $X$.
We will use this inequality in a moment to establish an estimate for the normal trace of $\alpha u$ and $\beta v$.
For this we will have to split the trace term $\norm{(u,v)}^2_{L^2(\Gamma)} \lesssim \int_\Gamma (\alpha u \cdot u + \beta v \cdot v) \dd \sigma$ appearing on the right-hand side of Corollary~\ref{cor:DivCurl2} into a normal and tangential component.

In the following, $u^\nu = (u\cdot\nu)\nu$ denotes the part of $u$ in normal direction, whereas $u^\tau = \nu\times(u\times\nu)$ is the tangential one. We can thus decompose $u$ as $u=u^\nu+u^\tau$.

\begin{Lemma}\label{lem:SplitNormTang}
	In the situation of Lemma~\ref{lem:DivCurl} above we have
	\[
		\int_s^t \int_\Gamma (\alpha u\cdot u + v\cdot \beta v)\dd \sigma \dd \tau \lesssim \int_s^t \int_\Gamma\left( \abs{\nu \cdot\alpha u}^2 +\abs{\nu \cdot\beta v}^2+ \abs{u^\tau}^2 + \abs{v^\tau}^2 \right) \dd \sigma \dd \tau\,.
	\]
\end{Lemma}

\begin{proof}
	We only treat $\tr(u\cdot \alpha u)$. The estimate for $v$ is shown analogously. Since the matrix operator $\alpha$ is assumed to be bounded and uniformly positive definite, we can compute
	\begin{align*}
		\int_s^t \int_\Gamma \alpha u\cdot u \dd \sigma \dd \tau & = \int_s^t \int_\Gamma \left(u^\nu\cdot\alpha u^\nu + 2u^\tau \alpha u^\nu+ u^\tau\cdot\alpha u^\tau\right)\dd \sigma \dd \tau \\
		 & \lesssim \int_s^t \int_\Gamma \left(\sqrt{\delta}\abs{u\cdot \nu} \frac{1}{\sqrt{\delta}}\abs{\nu \cdot\alpha u^\nu} + \abs{u^\tau}^2\right) \dd \sigma \dd \tau \\
		 & \lesssim \int_s^t \int_\Gamma \left(\delta \abs{u\cdot \nu}^2 + \frac{1}{\delta} \abs{\nu \cdot\alpha u^\nu}^2 + \abs{u^\tau}^2\right) \dd \sigma \dd \tau \\
		 & \lesssim \int_s^t \int_\Gamma \left(\delta \abs{u}^2 + \frac{1}{\delta}(\abs{\nu \cdot\alpha u}^2 +\abs{\nu \cdot\alpha u^\tau}^2) + \abs{u^\tau}^2\right) \dd \sigma \dd \tau \\
		 & \lesssim \int_s^t \int_\Gamma \left(\delta u\cdot\alpha u + \frac{1}{\delta}(\abs{\nu \cdot\alpha u}^2 +\abs{u^\tau}^2) + \abs{u^\tau}^2 \right) \dd \sigma \dd \tau\,.
	\end{align*}
	Fixing a sufficiently small $\delta>0$ we thus obtain the assertion.
\end{proof}

The next lemma allows us to control the normal trace of the curl terms. See also Chapter~2.3 in \cite{cessenat1996} for closely related results on the ``surface curl'' given by $\curl_\Gamma u_0 = \nu \cdot \curl u|_\Gamma$, where $u_0 = \tr_\tau u$.

\begin{Lemma}\label{Lem:surfaceCurl}
	For $f \in H^1(\Omega)$ we can estimate the normal surface curl by
	\[
		\norm{\nu \cdot( \curl f)}_{H^{-1}(\Gamma)} \lesssim \norm{\nu \times f}_{L^2(\Gamma)}\,.
	\]
\end{Lemma}

\begin{proof}
	Let $\phi \in H^1(\Gamma)$.
	This function can be extended in $H^{3/2}(\Omega)$ with $\norm{\phi}_{H^{3/2}(\Omega)} \lesssim \norm{\phi}_{H^1(\Gamma)}$.
	Let $f_n, \phi_m \in C_c^\infty(\overline\Omega)$ with $f_n \to f$ in $H^1(\Omega)$ and $\phi_m \to \phi$ in $H^{3/2}(\Omega)$.
	The divergence theorem and integration by parts yield
	\begin{align*}
		\int_\Gamma \nu \cdot(\curl f_n) \phi_m \dd \sigma & = (\curl f_n,\nabla\phi_m)_{L^2(\Omega)} + (\nabla \cdot (\curl f_n), \phi_m)_{L^2(\Omega)} \\
		 & = (f_n,\curl \nabla\phi_m)_{L^2(\Omega)} + \int_\Gamma (\nu \times f_n) \cdot \nabla \phi_m \dd \sigma \\
		 & \leq 0 + \norm{\nu \times f_n}_{L^2(\Gamma)}\norm{\tr_\tau(\nabla \phi_m)}_{L^2(\Gamma)} \\
		 &\lesssim \norm{\nu \times f_n}_{L^2(\Gamma)}\norm{\phi_m}_{H^1(\Gamma)}\,.
	\end{align*}
	Here the continuity of $\tr_\tau \nabla = \nabla_\Gamma \colon H^1(\Gamma) \to L^2(\Gamma)$ can be seen by writing the tangential gradient in terms of a parametrization of $\Gamma$, see for example Definition~2.3 and the subsequent remark in \cite{dziuk2013}.
	As $m \to \infty$ we obtain the estimate
	\[
		\norm{\nu \cdot( \curl f_n)}_{H^{-1}(\Gamma)} \lesssim \norm{\nu \times f_n}_{L^2(\Gamma)}\,.
	\]
	Letting $n$ tend to infinity then shows the claim.
\end{proof}

The main difficulty in the proof of Proposition \ref{prop:observability} is that one has to bound trace terms of the electric and magnetic field by the tangential trace of $\partial_t^k E$.
The boundary condition directly connects this trace with the tangential trace of the magnetic field.
We thus inspect the normal boundary terms in detail, where we gain small constants in front of the divergence terms, but obtain a (large) error in $H^{-1}(\Gamma)$.
In the next lemma, one can remove the operator $X$ on the right-hand side by continuity.
We keep it in the statement in view of Corollary~\ref{cor:inequality}.

\begin{Lemma}\label{Lem:normalTraces}
	Let $(u,v) \in G^1$ be a solution to the system \eqref{pde:2} and $\theta \in [0,\frac12)$.
	Then there exists a constant $\delta_0>0$ such that for every $\delta \in (0,\delta_0]$ and a constant $c_\delta$ depending on $\delta$ we get
	\begin{align}\label{eq:normalTraces}
		\int_{s}^{t} &\left(\norm{\nu \cdot \alpha u}_{L^2(\Gamma)}^2 + \norm{\nu \cdot \beta v}_{L^2(\Gamma)}^2\right) \dd \tau
		  \lesssim \int_{s}^{t} \Big(\delta\norm{u}_{H^{-\theta}(\div_\alpha)}^2 + \delta\norm{v}_{H^{-\theta}(\div_\beta)}^2  \notag\\
		  & \hphantom{<{}}+ c_\delta \norm{\nu \cdot (\alpha u, \beta v)}_{H^{-1}(\Gamma)}^2 + \norm{\nu\times u}_{L^2(\Gamma)}^2 + \norm{\omega}^2_{H^{1/2-\theta}(\Gamma)}\Big)\dd \tau \\
		  & \hphantom{<{}}+ \int_\R  \delta \norm{X( \xi, \zeta)}_{H^{-\theta}(\Omega)}^2  \dd \tau + \norm{(\nu \cdot \xi,\nu \cdot \zeta)}_{H^{-1}(\Gamma_\infty)}^2\notag
	\end{align}
	for $0 \leq s \leq t \leq T$.
\end{Lemma}

\begin{proof}
	We show the inequality only for $s=0$ and $t=T$ since general times can be treated analogously.

	As $\alpha u \in G^1([0,T])$, the normal trace belongs to $C([0,T],L^2(\Gamma))$. Extending $\alpha u$ by $0$ to $\R$, we obtain
	\begin{equation}\label{eq:partition}
		\begin{multlined}
		\int_0^{T} \norm{\nu \cdot \alpha u}_{L^2(\Gamma)}^2 \dd \tau
		= \norm{\nu \cdot \alpha u}_{L^2(\Gamma_\infty)}^2 \\
		\lesssim \norm{X(\nu \cdot \alpha u)}_{L^2(\Gamma_\infty)}^2
		+ \norm{(1-X)(\nu \cdot \alpha u)}_{L^2(\Gamma_\infty)}^2\,.
		\end{multlined}
	\end{equation}
	We first estimate the second term on the right-hand side.
	Lemma~\ref{lemma:eller:2.3} yields
	\[
		\norm{(1-X)(\nu \cdot \alpha u)}^2_{L^2(\Gamma_\infty)} \lesssim \norm{\nu\cdot\alpha u}^2_{H^{-1}(\Gamma_\infty)} + \norm{\partial_t (\nu \cdot \alpha u)}^2_{H^{-1}(\Gamma_\infty)}\,.
	\]
	We have $\nu \cdot \partial_t(\alpha u) = \nu\cdot(\curl v) + \nu \cdot \xi$ by \eqref{pde:2} and 
	\[
		\norm{\nu\cdot(\curl v)}_{H^{-1}(\Gamma_\infty)} \lesssim \norm{\nu \times v}_{L^2(\Gamma_\infty)} < \infty
	\]
	by Lemma~\ref{Lem:surfaceCurl}.
	Hence, the difference $\nu \cdot \xi$ belongs to $H^{-1}(\Gamma_T)$, and we can infer
	\begin{equation}\label{eq:proofOf4.9}
	\begin{multlined}
		\norm{(1-X)\nu \cdot (\alpha u)}^2_{L^2(\Gamma_\infty)}\\
		 \lesssim \int_0^{T} \left(\norm{\nu\cdot\alpha u}^2_{H^{-1}(\Gamma)} + \norm{\nu \times v}_{L^2(\Gamma)}^2 \right)\dd \tau + \norm{\nu \cdot \xi}_{H^{-1}(\Gamma_\infty)}^2
	\end{multlined}
\end{equation}
	By interpolation (see Theorem~2.7.7 in \cite{lions1972} for smooth $\Gamma$), the first term of the right-hand side of \eqref{eq:partition} can be bounded by
	\begin{align*}
		\norm{X (\nu \cdot \alpha u)}_{L^2(\Gamma_\infty)}^2
		 & \lesssim \int_\R \bigl(\hat\delta \norm{X (\nu \cdot \alpha u)}_{H^{1/2-\theta}(\Gamma)}^{2/(3-2\theta)}\frac{1}{\hat\delta}\norm{X (\nu \cdot \alpha u)}^{(1-2\theta)/(3-2\theta)}_{H^{-1}(\Gamma)} \bigr)^2\dd \tau\,.\\
		% \intertext{Now Young's product inequality with $p=3/2$ furnishes the estimate}
		 & \lesssim \int_\R \bigl(\hat\delta^{3-2\theta} \norm{X \nu \cdot \alpha u}_{H^{1/2-\theta}(\Gamma)}^{2}\\
		 & \qquad +\hat\delta^{-2(3-2\theta)/(1-2\theta)}\norm{X \nu \cdot \alpha u}^{2}_{H^{-1}(\Gamma)} \bigr)\dd \tau\,
	\end{align*}
	for an arbitrary $\delta \coloneqq \hat\delta^{3-2\theta} \in (0,1]$.
	Since $\Gamma \in C^5$ and $\alpha \in C^1$, the commutators $[X,\nu]$ and $[X,\alpha]$ are bounded from $L^2(\R,H^{-1/2-\theta})$ to $L^2(\R,H^{1/2-\theta})$, cf.\ Chapter~4 in \cite{taylor1991}.
	It follows that
	\begin{align*}
		 J &\coloneqq \int_\R (\delta \norm{X \nu \cdot \alpha u}_{H^{1/2-\theta}(\Gamma)}^{2}+c_\delta\norm{X \nu \cdot \alpha u}^{2}_{H^{-1}(\Gamma)} )\dd \tau \\
		 & \lesssim \int_\R \left(\delta \big(\norm{X u}^2_{H^{1/2-\theta}(\Gamma)} + \norm{u}^2_{H^{-1/2-\theta}(\Gamma)}\big) + c_\delta \norm{\nu \cdot \alpha u}_{H^{-1}(\Gamma)}^2\right) \dd \tau \\
		 & \lesssim \int_\R \left(\delta \big(\norm{X u}^2_{H^{1-\theta}(\Omega)} + \norm{u}^2_{H^{-1/2-\theta}(\Gamma)}\big) + c_\delta \norm{\nu \cdot \alpha u}_{H^{-1}(\Gamma)}^2\right) \dd \tau\,. \\
		 & \lesssim \int_\R \Big(\delta \big(\norm{u}_{H^{-\theta}(\div_\alpha)}^2 + \norm{v}_{H^{-\theta}(\div_\beta)}^2 + \norm{(u,v)}_{H^{-1/2-\theta}(\Gamma)}^2 + \norm{X(\xi,\zeta)}_{H^{-\theta}(\Omega)}^2 \\
		 &\qquad+ \norm{\omega}^2_{H^{1/2-\theta}(\Gamma)}\big)  + c_\delta \norm{\nu \cdot \alpha u}_{H^{-1}(\Gamma)}^2\Big) \dd \tau \,, \\
		\intertext{where we also used Corollary~\ref{cor:DivCurl2}. Lemma~\ref{lem:SplitNormTang} now allows us to split the trace term into normal and tangential parts, obtaining}
		J & \lesssim \delta \int_0^T\norm{\nu\cdot(\alpha u,\beta v)}_{L^2(\Gamma)}^2 \dd \tau+ \delta \int_\R \norm{X(\xi,\zeta)}_{H^{-\theta}(\Omega)}^2\dd\tau + \int_0^{T}\! \Big(\delta \big(\norm{u}_{H^{-\theta}(\div_\alpha)}^2\! \\
		 & \hphantom{<{}}+ \norm{v}_{H^{-\theta}(\div_\beta)}^2 \!+ \norm{\nu \times(u,v)}_{L^2(\Gamma)}^2 + \norm{\omega}^2_{H^{1/2-\theta}(\Gamma)}\big) + c_\delta \norm{\nu \cdot \alpha u}_{H^{-1}(\Gamma)}^2\Big) \dd \tau\,.
	\end{align*}
	We have the same estimate for $\nu \cdot \beta v$ in \eqref{eq:partition}.
	By choosing $\delta>0$ sufficiently small, we absorb the %$\norm{(\nu\cdot \alpha u, \nu\cdot \beta v)}_{L^2(\Gamma_\infty)}$
	normal trace terms by the left-hand side of \eqref{eq:partition}.
	We infer
	\begin{equation*}
		\begin{aligned}
		\int_0^{T}& \left(\norm{\nu \cdot \alpha u}_{L^2(\Gamma)}^2 + \norm{\nu \cdot \beta v}^2_{L^2(\Gamma)}\right) \dd \tau \lesssim \int_{0}^{T}\! \Big(\delta\norm{u}_{H^{-\theta}(\div_\alpha)}^2 \!+ \delta\norm{v}_{H^{-\theta}(\div_\beta)}^2 \\
		&+ c_\delta \norm{\nu \!\cdot\! (\alpha u, \beta v)}_{H^{-1}(\Gamma)}^2 + \norm{\nu\times (u,v)}_{L^2(\Gamma)}^2 + \delta \norm{\omega}^2_{H^{1/2-\theta}(\Gamma)} \Big)\dd \tau \\
		&+ \int_\R \delta \norm{X( \xi, \zeta)}_{H^{-\theta}(\Omega)}^2 \dd \tau
		+ \norm{(\nu \cdot \xi,\nu \cdot \zeta)}_{H^{-1}(\Gamma_\infty)}^2\,,
		\end{aligned}
	\end{equation*}
	using also \eqref{eq:proofOf4.9}.
	Finally, we note that the boundary condition in \eqref{pde:2} yields
	\begin{equation}\label{eq:tangentialEst}
		\begin{split}
			\int_{0}^{T} \norm{\nu \times v}_{L^2(\Gamma)}^2 \dd \tau \lesssim \int_{0}^{T} \left(\norm{\nu \times u}_{L^2(\Gamma)}^2 + \norm{\omega}_{L^2(\Gamma)}^2\right) \dd \tau\,,
		\end{split}
	\end{equation}
	from which the claim follows.
\end{proof}

We now rewrite the above result for the problem \eqref{pde:1} with $\theta = 0$ where $\xi=\partial_t\varphi$ and $\zeta = \partial_t \psi$.
Here we can apply Lemma~\ref{lemma:eller:2.2} to the inhomogeneities in \eqref{eq:normalTraces}.

\begin{Corollary}\label{cor:inequality}
	Let $\alpha$, $\beta$, $\varphi$, $\psi$, $\omega$, $u^{(0)}$, and $v^{(0)}$ satisfy \eqref{eq:regularityCoefficients} and \eqref{eq:regularityCoefficients2}.
	Then the solution $(u,v)\in C([0,T], L^2(\Omega))$ to \eqref{pde:1} possesses a trace in $L^2([0,T],L^2(\Gamma))$ satisfying
	\begin{equation}
		 \begin{aligned}\label{eq:47}
			 & \int_0^{T} \left(\norm{\nu \cdot \alpha u}_{L^2(\Gamma)}^2 + \norm{\nu \cdot \beta v}_{L^2(\Gamma)}^2\right)\dd \tau \\
			 & \qquad \lesssim \int_0^{T} \Big(\delta\norm{u}_{H(\div_\alpha)}^2 + \delta\norm{v}_{H(\div_\beta)}^2 + \norm{(\varphi, \psi)}_{H^1(\Omega)}^2 \\
			 & \qquad\qquad + c_\delta\norm{\nu \cdot (\alpha u, \beta v)}_{H^{-1}(\Gamma)}^2 + \norm{\nu\times u}_{L^2(\Gamma)}^2 + \norm{\omega}^2_{H^{1/2}(\Gamma)} \Big)\dd \tau \\
			 & \hphantom{<{}} \qquad + \norm{(\nu \cdot \partial_t \varphi,\nu \cdot \partial_t \psi)}_{H^{-1}(\Gamma_\infty)}^2\,.
		 \end{aligned}	 		
	\end{equation}
\end{Corollary}

\begin{proof}
	For $(u,v) \in G^1$ the assertion follows from Lemma~\ref{Lem:normalTraces} with $\xi = \partial_t \varphi$ and $\zeta = \partial_t \psi$, where we apply Lemma~\ref{lemma:eller:2.2} to the term with $X$ on the right-hand side.
	It remains to remove the assumption $(u,v) \in G^1$.
	This can be done as in the proof of Lemma~3.4 in \cite[]{pokojovy2020}.
	One approximates the initial data in $H(\div_{\alpha(0)}) \times H(\div_{\beta(0)})$ by $(u^{(0)}_n, v^{(0)}_n) \in H^1(\Omega)$ satisfying the compatibility condition
	$v^{(0)}_n \times \nu + (\gamma(u^{(0)}_n\times\nu))\nu  = 0$,
	$(\varphi, \psi)$ in $G^1$ by $(\varphi_n, \psi_n) \in G^2$, and 
	$\omega$ in $L^2([0,T], H^{1/2}(\Gamma))$ by $\omega_n \in H^1([0,T] \times \Gamma)$ with $\omega \cdot \nu = 0$ and $\omega_n(0)=0$. 
	By Theorem 1.3 and Remark 2.1 in \cite[]{cagnol2011} there are solutions $(u_n,v_n) \in G^1$.
	%with $(\tr u_n, \tr v_n)$ in $H^1([0,T]\times\Gamma)$.
	According to Proposition~1.1 in \cite[]{cagnol2011} these solutions converge to $(u,v)$ in $C([0,T],L^2(\Omega))$ and  $(\tr_t u_n, \tr_t v_n)$ to $(\tr_t u, \tr_t v)$ in $L^2([0,T] \times \Omega)$.
	Moreover, we have the limits $\div(\alpha u_n) = \div(\varphi_n) \to \div(\alpha u)$ and $\div(\beta v_n) = \div(\psi_n) \to \div (\beta v)$ in $C([0,T], L^2(\Omega))$.
	As a result the traces $\tr_n (\alpha u_n)$ and $\tr_n(\beta v_n)$ converge in $C([0,T], H^{-1/2}(\Gamma))$.
	The assertion for $(u,v)$ now follows from the claim for $(u_n,v_n)$.
\end{proof}

\begin{Remark}\label{rem:est}
	In Corollary~\ref{cor:inequality} we can insert $\div(\alpha u) = \div(\varphi)$ and $\div(\beta v) = \div(\psi)$
	and estimate
	\[
		\norm{(\nu \cdot \partial_t \varphi,\nu \cdot \partial_t \psi)}_{H^{-1}(\Gamma_\infty)}^2 \lesssim \norm{(\nu \cdot \varphi,\nu \cdot \psi)}_{L^2([0,T], L^2(\Gamma))}^2\lesssim \norm{(\varphi,\psi)}_{L^2([0,T], H^1(\Omega))}^2\,.
	\]
	Moreover, we can bound $c_\delta\norm{\nu \cdot (\alpha u, \beta v)}_{H^{-1}(\Gamma)}^2$ by
	\[
		c_\delta\big(\norm{u}_{H(\div_\alpha)}^2 + \norm{v}_{H(\div_\beta)}^2\big)\,.
	\]
\end{Remark}

The estimate~\eqref{eq:47} has time independent constants and can therefore be used below to show Proposition~\ref{prop:observability} and thus Theorem~\ref{thm:mainThm}.
Concluding this section, we also prove a variant only having norms of data on the right-hand side but with time dependent constants.
In this way we give a different proof of Theorem~1.2 in \cite{cagnol2011}.
We will revisit this trace estimate in Section~\ref{sec:6} in the linear, autonomous and homogeneous case and prove a version with time independent constants there.

In view of Remark~\ref{rem:est} we have to remove the $L^2$-norm of $u$ and $v$ and the tangential trace term on the right-hand side of \eqref{eq:normalTraces}. We control these terms by means of the energy estimate for \eqref{pde:1}.
It was shown in Lemma~3.2 of \cite[]{pokojovy2020} for $\omega=0$.
An obvious modification of this proof gives the equality
\begin{equation}\label{eq:energy_inequality}
	\begin{aligned}
	&\norm{\alpha(t)^{1/2} u(t)}_{L^2}^2  + \norm{\beta(t)^{1/2}v(t)}_{L^2}^2 + 2 \int_0^t\norm{\gamma^{1/2} \tr_t u(\tau)}_{L^2(\Gamma)}^2 \dd\tau \\
	 & \quad =\norm{\alpha(0)^{1/2} u(0)}_{L^2}^2+\norm{\beta(0)^{1/2} v(0)}_{L^2}^2+2\int_{0}^{t}\int_{\Omega}(u\cdot \partial_t\varphi+v\cdot \partial_t\psi)\,\dd x \dd \tau \\
	 & \qquad +\int_0^t\int_\Gamma \tr_t u \cdot \omega \dd \sigma \dd \tau + \int_{0}^{t}\int_{\Omega} \left(u\cdot(\partial_t\alpha)u+v\cdot(\partial_t\beta)v\right) \dd x \dd \tau \,,
\end{aligned}
\end{equation}
for $0 \leq t \leq T$.
To remove the trace term on the right-hand side, we estimate
\begin{equation}
	\abs{\int_0^t\int_\Gamma \tr_t u \cdot \omega \dd \sigma \dd \tau} \leq \delta \int_0^t \norm{\tr_t u}_{L^2(\Gamma)}^2 \dd \tau + c_\delta \int_0^t \norm{\omega}^2_{L^2(\Gamma)} \dd \tau\,.
\end{equation}
Since $\gamma  \geq \eta$, setting $\delta = \eta$ we arrive at
	\begin{multline}
		\norm{(u(t), v(t))}^2_{L^2} + \int_0^t\norm{\tr_t u(\tau)}^2_{L^2(\Gamma)} \dd \tau \leq c \norm{(u(0), v(0))}^2_{L^2}\\ 
		+ \overline{c} \int_0^t \norm{(u,v)}^2_{L^2} \dd \tau 
		 + c \int_0^t \left(\norm{(\partial_t \varphi, \partial_t \psi)}^2_{L^2} + \norm{\omega}^2_{L^2(\Gamma)}\right) \dd \tau
	\end{multline}
with $\overline{c} \coloneqq \tilde c + \max\{\norm{\partial_t \alpha}_\infty , \norm{\partial_t \beta}_\infty\}$ where $\tilde c=0$ if $\partial_t \varphi = \partial_t \psi =0$ and $\tilde c=1$ otherwise.
Gronwall's inequality now implies
\begin{multline}
		\norm{(u(t), v(t))}^2_{L^2} + \int_0^t\norm{\tr_t u(\tau)}^2_{L^2(\Gamma)} \dd \tau \\
		\lesssim \bigg( \norm{(u(0),v(0))}^2_{L^2}  + \int_0^t \left(\norm{(\partial_t \varphi, \partial_t \psi)}^2_{L^2} + \norm{\omega}^2_{L^2(\Gamma)}\right) \dd \tau \bigg) \e^{\overline{c}t}
\end{multline}
Using also $\div(\alpha u) = \div(\varphi)$ and $\div(\beta v) = \div(\psi)$, from Corollary~\ref{cor:inequality} and Remark~\ref{rem:est} we derive a trace estimate for the inhomogeneous, linear Maxwell system \eqref{pde:1}. It can directly be applied to \eqref{eq:maxwell}--\eqref{eq:maxwell4}.

\begin{Proposition}
	\label{prop:gronwall}
	Assume that conditions \eqref{eq:regularityCoefficients} and \eqref{eq:regularityCoefficients2} hold and let $(u,v) \in C([0,T], L^2(\Omega))$ solve \eqref{pde:1}. We obtain
	\begin{equation*}
		\begin{split}
		\int_0^t \norm{\tr(u,v)}^2_{L^2(\Gamma)} \dd \tau
		&\lesssim_T \norm{(u^{(0)}, v^{(0)})}^2_{L^2(\Omega)} \\
		&\hphantom{{}=}
		+ \int_0^t 
		\Big(\norm{(\partial_t\varphi, \partial_t \psi)}^2_{L^2(\Omega)}
		+ \norm{(\varphi, \psi)}^2_{H^1(\Omega)}
		+ \norm{\omega}^2_{L^2(\Gamma)}\Big) 
		\dd \tau			
	\end{split}
	\end{equation*}
	for $0 \leq t \leq T$.
\end{Proposition}

As in Theorem~1.2 of \cite[]{cagnol2011} the constant depends on $T$. Thus, we will only use Corollary~\ref{cor:inequality} in the following.

\section{Proof of Proposition \ref{prop:observability}}
\label{sec:5}

We again consider system \eqref{pde:1} assuming \eqref{eq:regularityCoefficients}, \eqref{eq:regularityCoefficients2}, and
%With Corollary \ref{cor:inequality} we can follow the proof of Lemma~3.4 in \cite{pokojovy2020} where $\omega = 0$. By a regularization argument it was shown that it suffices to consider $(u,v) \in G^1([0,T])$ and $\varphi, \psi \in G^2$. From the proof we furthermore get the following estimate under the assumption that $\alpha$ and $\beta$ satisfy
\begin{equation}\label{eq:nontrapping}
	\alpha + (m \cdot \nabla)\alpha \geq  \tilde{\eta} \alpha, \quad
	\beta + (m \cdot \nabla)\beta \geq \tilde{\eta} \beta
\end{equation}
for some $\tilde{\eta}>0$.
The proof of Lemma~3.4 in \cite{pokojovy2020} then yields the following result.
In \cite{pokojovy2020} it was assumed that $\Omega$ is strictly starshaped giving control on the full trace of $(u,v)$ which now appears on the right-hand side of \eqref{eq:pokojovyEst}.
Moreover, the boundary inhomogeneity $\omega$ can simply be included in the calculations of \cite{pokojovy2020}.
\begin{Lemma} \label{lem:InnerEnergyEst}
	Let \eqref{eq:regularityCoefficients}, \eqref{eq:regularityCoefficients2} and \eqref{eq:nontrapping} hold.
	Then the weak solutions $(u,v) \in C([0,T], L^2(\Omega,\R^6))$ of system \eqref{pde:1} fulfill
	\begin{equation}\label{eq:pokojovyEst}
	\begin{split}
		 & \int_s^t \int_\Omega (\alpha u \cdot u + \beta v \cdot v) \dd \sigma \dd \tau \\
		 & \quad \lesssim \int_s^t \int_\Gamma \abs{u \times \nu}^2 \dd \sigma \dd \tau + \norm{(u(s),v(s))}^2_{L^2}+ \norm{(u(t),v(t))}^2_{L^2} \\
		 & \qquad + \int_s^t \int_{\Omega} \left(\abs{\partial_t\varphi} \abs{v} + \abs{\partial_t \psi} \abs{u} + \abs{\div\varphi} \abs{u} + \abs{\div\psi} \abs{v}\right) \dd x \dd \tau\\
		 & \qquad + \int_s^t \int_\Gamma \tr(\alpha u \cdot u + \beta v \cdot v) \dd \sigma \dd \tau + \int_s^t \int_\Gamma \abs{\omega \cdot \tr v} \dd \sigma \dd \tau
	\end{split}
\end{equation}
	for $0\leq s \leq t \leq T$.
	Note that the last term is bounded by $c \int_s^t \int_\Gamma\tr(\beta v \cdot v) \dd \sigma \dd \tau + c \int_s^t \norm{\omega}^2_{L^2(\Gamma)} \dd \tau$ and that the trace terms are finite by Proposition~\ref{prop:gronwall}.
\end{Lemma}

We still have to estimate the integral $\int_s^t \int_\Gamma \tr(\alpha u \cdot u + \beta v \cdot v )\dd \sigma \dd \tau$ uniformly in $T$.
By Lemma~\ref{lem:SplitNormTang} it can be split into 
	\[
		\int_s^t \int_\Gamma \tr(\alpha u \cdot u + \beta v \cdot v )\dd \sigma \dd \tau \lesssim \int_s^t \int_\Gamma \left(\abs{\nu \cdot\alpha u}^2 + \abs{u^\tau}^2 + \abs{\nu \cdot\beta v}^2 + \abs{v^\tau}^2\right) \dd \sigma \dd \tau\,.
	\]
	The tangential trace of $v$ can be controlled through the boundary condition by
	\begin{equation*}
		\int_{s}^{t} \norm{v^\tau}^2_{L^2(\Gamma)} \dd \tau \lesssim \int_{s}^{t}\left(\norm{\nu \times u}^2_{L^2(\Gamma)} + \norm{\omega}^2_{L^2(\Gamma)}\right) \dd \tau\,.
	\end{equation*}
Using Corollary \ref{cor:inequality} and Remark~\ref{rem:est}, we derive
	\begin{multline*}
		\int_{s}^{t} \int_\Gamma \left(\abs{\nu \cdot \alpha u}^2+\abs{\nu \cdot \beta v}^2\right) \dd \sigma \dd \tau
		  \lesssim \int_s^t \Big(\delta\norm{u}_{L^2(\Omega)}^2 + \delta\norm{v}_{L^2(\Omega)}^2 + \norm{(\varphi,\psi)}_{H^1(\Omega)}^2 \\
		 \qquad + \norm{\nu\times u}_{L^2(\Gamma)}^2+ \norm{\omega}^2_{H^{1/2}(\Gamma)} + c_{\delta}\norm{\nu\cdot(\alpha u, \beta v)}^2_{H^{-1}(\Gamma)} \Big) \dd \tau\,.
	\end{multline*}
	Fixing a sufficiently small $\delta>0$, the norms in $L^2$ above can be absorbed by the left-hand side of \eqref{eq:pokojovyEst}. We have shown the following estimate.

	\begin{Lemma}\label{lem:conclusion}
		Let \eqref{eq:regularityCoefficients}, \eqref{eq:regularityCoefficients2} and \eqref{eq:nontrapping} be true.
		Then the weak solutions $(u,v) \in C([0,T], L^2(\Omega,\R^6))$ of system \eqref{pde:1} fulfill
		\begin{align}\label{eq:conclusion}
		\begin{split}
		 \int_s^t& \int_\Omega (\alpha u \cdot u + \beta v \cdot v )\dd x \dd \tau + \int_{s}^{t} \int_\Gamma \tr(\alpha u \cdot u + \beta v \cdot v ) \\
		 &\lesssim \int_s^t \int_\Gamma \abs{u \times \nu}^2 \dd \sigma \dd \tau + \norm{(u(s),v(s))}^2_{L^2}+ \norm{(u(t),v(t))}^2_{L^2} \\
		 & \qquad + \int_s^t \int_{\Omega}\bigl( \abs{\partial_t\varphi} \abs{v} + \abs{\partial_t \psi} \abs{u} + \abs{\div\varphi} \abs{u} + \abs{\div\psi} \abs{v}\bigr) \dd x \dd \tau \\
		 & \qquad+ \int_s^t \bigl( \norm{(\varphi , \psi)}_{H^1}^2 + \norm{\omega}_{H^{1/2}(\Gamma)}^2 + \norm{\nu \cdot (\alpha u, \beta v)}_{H^{-1}(\Gamma)}^2 \bigr)\dd \tau
		\end{split}
	\end{align}
		for $0\leq s \leq t \leq T$.
	\end{Lemma}

One could estimate the terms in the second line of the right-hand side by 
\[
	\delta \int_s^t \int_\Omega \bigl(\abs{u}^2 + \abs{v}^2\bigr) \dd x \dd \tau + c_\delta \int_s^t \int_\Omega \bigl(\abs{\partial_t\varphi}^2 + \abs{\partial_t\psi}^2 + \abs{\div \varphi}^2 + \abs{\div \phi}^2 \bigr)\dd x \dd \tau \,,
\]
and absorb the first summand for small $\delta>0$ by the left-hand side.
In the following this is not needed, since these summands can be put into small error terms.

In order to obtain Proposition \ref{prop:observability} it remains to control the $H^{-1}$-norm of the normal traces.
To this end, we return to our original problem.
Here it does not help to bound them by the norm of $u$ and $v$ in $H(\div_\alpha)$ and $H(\div_\beta)$, since the time integrals of $L^2$-norms of $u$ and $v$ cannot be absorbed by the left-hand side.

We crucially exploit that the time derivatives of solutions still solve a Maxwell system.
We begin with $\left(\partial_t^k E, \partial_t^k H\right)$ for $k\geq1$ and obtain estimates with dissipation terms $d_{k-1}$ on the right-hand side.
In turn the estimate for $(\partial_t E, \partial_t H)$ from Lemma~\ref{lem:RandtermH-1k123} allows us to treat $(E,H)$ in Lemma~\ref{lem:RandtermH-1k0}, where we also use the invertibility of curl in our setting to be proved in Theorem~\ref{thm:div-curl-est}~b).
To this purpose we assume again that $\R^3\setminus\Omega$ is connected.

\begin{Lemma} \label{lem:RandtermH-1k123}
	Let the assumptions of Theorem~\ref{thm:mainThm} hold.
	For the solution $(E,H) \in G^3$ of \eqref{eq:maxwell}--\eqref{eq:maxwell4} we can estimate the normal trace by
	\begin{align*}
		\int_s^t\norm{\nu \cdot (\widehat{\epsilon}_k \partial_t^k E, \widehat{\mu}_k \partial_t^k H)}^2_{H^{-1}(\Gamma)} \dd\tau
		\lesssim \int_s^t \left(d(\tau) + z(\tau)^2\right) \dd \tau
	\end{align*}
	for $0 \leq s \leq t < T_*$ and $k \in \{1,2,3\}$.
\end{Lemma}

\begin{proof}
	Let $k \in \{1,2,3\}$. Due to \eqref{eq:derivative} we can rewrite $\widehat{\epsilon}_k\partial_t^k E$ as
	\begin{align*}
		\widehat{\epsilon}_k\partial_t^k E & = \partial_t^k(\epsilon E) - f_k = \partial_t(\widehat{\epsilon}_{k-1}\partial_t^{k-1} E + f_{k-1}) - f_k\,.
		\intertext{The Maxwell system \eqref{eq:inhomogenous_maxwell} then leads to}
		\widehat{\epsilon}_k\partial_t^k E &= \curl{\partial_t^{k-1}H} - \partial_t f_{k-1} + \partial_t f_{k-1} -f_k = \curl{\partial_t^{k-1}H}-f_k
	\end{align*}
	Let $\Phi \in H^1(\Gamma)$. The equations above and Lemma~\ref{Lem:surfaceCurl} yield
	\begin{align*}
		\langle \nu \cdot \widehat{\epsilon}_k\partial_t^k E, \Phi \rangle_{H^{-1/2}(\Gamma)} & = \langle \nu \cdot (\curl \partial_t^{k-1}H), \Phi \rangle_{H^{-1}(\Gamma)} - \langle \nu \cdot f_k , \Phi \rangle_{L^2(\Gamma)} \\
		 & \lesssim \norm{\nu \times \partial_t^{k-1}H}_{L^2(\Gamma)}\norm{\Phi}_{H^1(\Gamma)} + \norm{f_k}_{L^2(\Gamma)}\norm{\Phi}_{L^2(\Gamma)}\,.
	\end{align*}
	Dividing by $\norm{\Phi}_{H^1(\Gamma)}$, we infer
	\begin{align*}
		\norm{\nu \cdot \widehat{\epsilon}_k \partial_t^k E}_{H^{-1}(\Gamma)} &\lesssim \norm{\nu \times \partial_t^{k-1}H}_{L^2(\Gamma)} + \norm{f_k}_{H^1(\Omega)} \\
		&\lesssim \norm{\nu \times \partial_t^{k-1}E}_{L^2(\Gamma)} + \norm{h_{k-1}}_{L^2(\Gamma)} + \norm{f_k}_{H^1(\Omega)}\,,
	\end{align*}
	where we also used the boundary condition \eqref{eq:inhomogenous_maxwell_boundary}. Analogously we treat $\widehat{\mu}_k\partial_t^k H$.
	Estimate \eqref{eq:commutator-estimates} now implies the assertion.
\end{proof}

\begin{Lemma}
	\label{lem:RandtermH-1k0}
	Let the assumptions of Theorem~\ref{thm:mainThm} hold.
	For the solution $(E,H) \in G^3$ of \eqref{eq:maxwell}--\eqref{eq:maxwell4} we estimate
	\begin{equation*}
		\int_s^t\norm{\nu \cdot (\widehat{\epsilon}_0 E, \widehat{\mu}_0 H)}^2_{H^{-1}(\Gamma)} \dd \tau \lesssim \int_s^t d(\tau) \dd \tau + (e_1(t) + e_1(s)) + \int_s^t z^2(\tau) \dd \tau
	\end{equation*}
	for $0 \leq s \leq t < T_*$.
\end{Lemma}

\begin{proof}
	The usual normal trace estimate, \eqref{eq:maxwell2}, Theorem~\ref{thm:div-curl-est}~b), \eqref{eq:derivative} and \eqref{eq:inhomogenous_maxwell} yield
	\begin{align*}
		\int_s^t\norm{\nu \cdot (\widehat{\epsilon}_0 E, \widehat{\mu}_0 H)}^2_{H^{-1}(\Gamma)} \dd \tau & \lesssim \int_s^t\norm{(\widehat{\epsilon}_0 E, \widehat{\mu}_0 H)}^2_{H(\div)} \dd \tau \\
		&= \int_s^t\norm{(\widehat{\epsilon}_0 E, \widehat{\mu}_0 H)}^2_{L^2(\Omega)} \dd \tau\\
		%& \lesssim \int_s^t\norm{(E,H)}^2_{H^1(\Omega)} \dd \tau \\
		 & \lesssim \int_s^t\norm{(\curl E, \curl H)}^2_{L^2(\Omega)} \dd \tau \\
		%& = \int_s^t\norm{\partial_t(\widehat{\epsilon}_0 E, \widehat{\mu}_0 H)}^2_{L^2(\Omega)} \dd \tau \\
		 & = \int_s^t\norm{(\widehat{\epsilon}_1 \partial_t E, \widehat{\mu}_1 \partial_t H)}^2_{L^2(\Omega)} \dd \tau \\
		 & \lesssim \int_s^t \left(\partial_t E \cdot \widehat{\epsilon}_1 \partial_t E + \partial_t H \cdot \widehat{\mu}_1 \partial_t H \right) \dd \tau\,.
	\end{align*}
	Using Lemma~\ref{lem:conclusion} with $u= \partial_t E$, $v = \partial_t H$, $\alpha = \widehat{\epsilon}_1$, $\beta = \widehat{\mu}_1$ and $\varphi = \psi = \omega = 0$, % as well as \eqref{eq:derivative},
	we then obtain
	\begin{align*}
		 \int_s^t\norm{\nu \cdot (\widehat{\epsilon}_0 E, \widehat{\mu}_0 H)}^2_{H^{-1}(\Gamma)} \dd \tau 
		 & \lesssim \int_s^t d(\tau) \dd \tau + (e_1(t) + e_1(s)) \\%+ \int_s^t z^{3/2}(\tau)\dd \tau \\
		  &\quad + \int_s^t \norm{\nu\cdot(\widehat{\epsilon}_1 \partial_t E, \widehat{\mu}_1 \partial_t H)}^2_{H^{-1}(\Gamma)} \dd \tau\,.
	\end{align*}
	By the previous lemma the last term can be estimated by $\int_s^t \left(d(\tau) + z(\tau)^2\right)\dd \tau$ and the claim follows.
\end{proof} 

Combining the lemmas above and \eqref{eq:commutator-estimates}, we have shown Proposition~\ref{prop:observability}.

\section{The autonomous linear case} \label{sec:time-indep-bnds}
\label{sec:6}

Next we will improve the trace estimate from \cite{cagnol2011} or Proposition~\ref{prop:gronwall} in the linear, autonomous and homogeneous case, where we obtain $T$-independent constants.
Moreover we show exponential stability of this problem.
So we consider
\begin{align} \label{pde:lin}
	\notag	\partial_t (\epsilon(x) E)                           & = \curl H\,,                                  \\
	\notag	\partial_t (\mu(x) H)                                & = - \curl E\,,\qquad &&t\geq 0, x \in \Omega\,, \\
	\div(\varepsilon E) = 0                                          & = \div(\mu H)                                      \\
	\notag H \times \nu + \big(\lambda(x)(E\times\nu)\big)\times \nu & = 0\,,   \qquad &&t\geq 0, x \in \Gamma\,,                                          \\
	\notag E(0) = E^{(0)}\,, \quad H(0)                              & =H^{(0)}\,,\qquad &&x \in \Omega\,.
\end{align}
As in \eqref{eq:regularityCoefficients}, \eqref{eq:regularityCoefficients2} and \eqref{eq:nontrapping} we assume that
\begin{equation}\label{eq:cond-on-coefficients}
	\epsilon, \mu \in C^1\left(\overline{\Omega}, \Rsym\right),\ \lambda \in C^1_\tau\left(\Gamma, \Rsym\right), \quad \text{with } \epsilon, \mu, \lambda \geq \eta\,,
\end{equation}
for some $\eta>0$ and that
\begin{equation}\label{eq:cond-on-initial-values}
	E^{(0)}, H^{(0)} \in (L^2(\Omega))^3 \quad\text{satisfy}\quad \div(\varepsilon E^{(0)}) = \div(\mu H^{(0)})=0\,.
\end{equation}
Furthermore, we require
\begin{equation}\label{eq:non-trapping}
	\epsilon + (m \cdot \nabla)\epsilon \geq  \tilde{\eta} \epsilon, \quad
	\mu + (m \cdot \nabla)\mu \geq \tilde{\eta} \mu
\end{equation}
for some $\tilde \eta>0$.
We assume in this section that the complement
\begin{equation}\label{eq:complement-connected}
	\R^3\setminus\Omega \text{ is connected.}
\end{equation}

On $X \coloneqq \{ (E,H) \in L^2(\Omega)^6 \mid \div(\epsilon E) = \div(\mu H) = 0\}$ we define the operator
\[
	A \coloneqq \begin{pmatrix}
		0 & \varepsilon^{-1} \curl \\ -\mu^{-1} \curl & 0
	\end{pmatrix}
\]
with
\[
	D(A) \coloneqq \{ (E,H) \in X  \mid \curl E, \curl H \in L^2(\Omega), H \times \nu + (\lambda(E\times\nu))\times \nu = 0\}\,.
\]
Then \eqref{pde:lin} reduces to the evolution equation $\partial_t \begin{psmallmatrix}
		E \\H
	\end{psmallmatrix} = A\begin{psmallmatrix}
		E \\H
	\end{psmallmatrix}$.
As in Lemma~7.2.2.1 of \cite{lagnese2004} one sees that $A$ is maximally dissipative and thus generates a contraction semigroup $(S(t))_{t \geq 0}$.
We recall the energy estimate from equation~(7.2.2.12) of \cite{lagnese2004} or Lemma~3.2 of \cite{pokojovy2020}, namely
\begin{equation}\label{eq:energy-est-lin}
	\begin{aligned}
		 & \norm{(\epsilon^{1/2}E,\mu^{1/2}H)(s)}_{L^2}^2 - \int_s^t \norm{\lambda \tr_t E \cdot \tr_t E }_{L^2(\Gamma)}^2 \dd \tau \\
		 & \quad= \norm{(\epsilon^{1/2}E,\mu^{1/2}H)(t)}_{L^2}^2\,.
	\end{aligned}
\end{equation}
for $t\geq s \geq 0$.
Note that in contrast to \cite{cagnol2011} we do not consider time dependent coefficients $\epsilon$ and $\mu$, since this would for example introduce the extra term $\int_s^t \int_\Omega E \cdot (\partial_t \epsilon) E + H \cdot (\partial_t \mu) H$ on the left-hand side of \eqref{eq:energy-est-lin}.
%sign-error in \cite{pokojovy2020} there.

In combination with Lemma~\ref{lem:InnerEnergyEst}, the energy estimate above leads to the following preliminary observability inequality.
\begin{Lemma}\label{lem:observability}
	Let \eqref{eq:cond-on-coefficients}, \eqref{eq:cond-on-initial-values}, \eqref{eq:non-trapping} and \eqref{eq:complement-connected} be true, $(E^{(0)}, H^{(0)}) \in X$, and let $(E,H)\in C([0,\infty), X)$ solve \eqref{pde:lin}.
	Then there is a number $\hat T >0$ not depending on $(E^{(0)}, H^{(0)})$ such that for $T \geq \hat T$ we have the observability estimate
	\begin{equation}\label{eq:proof-time-independent-bnd-obs-1}
		\norm{(E,H)(0)}_{L^2}^2 \lesssim \int_0^T \int_\Gamma\tr(\epsilon E \cdot E + \mu H \cdot H) \dd \sigma \dd \tau \,.
	\end{equation}
\end{Lemma}

Hence, if $E$ and $H$ vanish on the boundary, the initial values $(E,H)(0)$ were zero and so, since we are in the homogeneous case, the fields $E$ and $H$ vanish as well.

\begin{proof}
	Lemma~\ref{lem:InnerEnergyEst} yields the observability-like estimate
	\begin{equation}\label{eq:Lem-5.1}
		\begin{aligned}
			 & \int_s^t\int_\Omega (\epsilon E \cdot E + \mu H \cdot H) \dd x \dd \tau
			\lesssim \norm{(E,H)(s)}_{L^2}^2 + \norm{(E,H)(t)}^2_{L^2}               \\
			 & \quad
			+ \int_s^t\int_\Gamma \tr(\epsilon E \cdot E + \mu H \cdot H) \dd \sigma \dd \tau %+\int_s^t \abs{E \times \nu}^2 \dd \tau
		\end{aligned}
	\end{equation}
	for $s,t \geq 0$.
	We stress that the implicit constant does not depend on t or s.

	Note that $\abs{\tr E}^2 \geq \abs{\tr_\tau E}^2 = \abs{\tr_t E}^2$ by orthogonality. The boundary condition in \eqref{pde:lin} as well as inequalities \eqref{eq:energy-est-lin} and \eqref{eq:Lem-5.1} then lead to
	\begin{align*}
		 & T\Big(\norm{(E,H)(0)}_{L^2}^2  - \int_0^T\int_\Gamma\tr(\epsilon E \cdot E + \mu H \cdot H) \dd \sigma\dd \tau\Big)\\
		 &\quad \lesssim T\Big(\norm{(E,H)(0)}_{L^2}^2  - \int_0^T\int_\Gamma\abs{\tr_t E}^2 \dd \sigma\dd \tau\Big)\\
		 & \quad \lesssim T \norm{(E,H)(T)}_{L^2}^2
		\lesssim \int_0^T \norm{(E,H)(\tau)}_{L^2}^2 \dd \tau\\
		 & \quad \lesssim \norm{(E,H)(0)}^2_{L^2}+ \norm{(E,H)(T)}^2_{L^2} + \int_0^T\int_\Gamma\tr(\epsilon E \cdot E + \mu H \cdot H) \dd \sigma\dd \tau \\
		 & \quad \lesssim \norm{(E,H)(0)}^2_{L^2} + \int_0^T\int_\Gamma\tr(\epsilon E \cdot E + \mu H \cdot H) \dd \sigma\dd \tau \,,
	\end{align*}
	cf.\ \eqref{eq:obs-est-linear}. The desired estimate follows for sufficiently large $T$.
\end{proof}

We now bound the normal traces of $E$ and $H$ solely through the tangential trace of the electric field, where the constants do not depend on $T$.
This improves Lemma~\ref{Lem:normalTraces} considerably in the autonomous, homogeneous case.
We employ a compactness argument which was originally used in the proof of global decay rates for the wave equation, see Lemma~2.2 in \cite{lasiecka1992-1} (see also \cite{eller2007}).
 Moreover, we use uniqueness properties both on the dynamical and stationary level which rely on the observability estimates in Lemma~\ref{lem:observability} and on the div-curl bound in Theorem~\ref{thm:div-curl-est}~b).

\begin{Theorem}\label{thm:normal-trace-vs-tang-trace}
	Let \eqref{eq:cond-on-coefficients}, \eqref{eq:cond-on-initial-values}, \eqref{eq:non-trapping}, and \eqref{eq:complement-connected} hold. Then the solution $(E,H) \in C([0,\infty), L^2(\Omega)^6)$ of \eqref{pde:lin} satisfies
	\begin{equation}
		\int_{s}^{t} \left(\norm{\nu \cdot \epsilon E}_{L^2(\Gamma)}^2 + \norm{\nu \cdot \mu H}_{L^2(\Gamma)}^2\right) \dd \tau
		\lesssim \int_{s}^{t} \norm{\nu\times E}_{L^2(\Gamma)}^2 \dd \tau
	\end{equation}
	for $0\leq s \leq s + \hat T \leq t$, where $\hat T$ is given by Lemma~\ref{lem:observability}.
\end{Theorem}

\begin{proof}
	1) Let $(E^{(0)}, H^{(0)}) \in X$ and $(E,H)$ be the corresponding solution of \eqref{pde:lin}. Take $T \geq \hat T$ with $\hat T$ from Lemma~\ref{lem:observability}.
	By Proposition~\ref{prop:gronwall} it suffices to choose $(E^{(0)}, H^{(0)}) \in D(A)$.

	We first show that for a fixed $T\geq\hat T$ we can estimate
	\begin{equation} \label{eq:proof-exp-stab:normal-trace}
		\int_0^T \norm{\nu \cdot (\epsilon E, \mu H)}_{L^2(\Gamma)}^2 \dd \tau
		\lesssim_T \int_0^T \norm{\nu\times E}_{L^2(\Gamma)}^2 \dd \tau\,.
	\end{equation}
	% (Of course we only need this estimate for the $H^{-1}$-norm in the previous inequality \eqref{eq:proof-exp-stab}.)
	Let $\theta \in (0,\frac12)$. In Lemma~\ref{Lem:normalTraces} we have already shown
	\begin{align}\label{eq:proof-time-independent-bnd}
		\notag & \int_0^T \left(\norm{\nu \cdot \epsilon E}_{L^2(\Gamma)}^2 + \norm{\nu \cdot \mu H}_{L^2(\Gamma)}^2\right) \dd \tau \\
		       & \quad \lesssim \int_0^T  \big(\delta\norm{(E,H)}_{H^{-\theta}}^2
		+  c\norm{\nu\times E}_{L^2(\Gamma)}^2 + c_\delta \norm{\nu \cdot (\epsilon E, \mu H)}_{H^{-1}(\Gamma)}^2\big)\dd \tau\,,
	\end{align}
	where we have fixed $T\geq \hat T$.

	2.a) Let us assume that \eqref{eq:proof-exp-stab:normal-trace} does not hold.
	Thus, there exist solutions $(\tilde E_n, \tilde H_n)$ satisfying the divergence condition in \eqref{pde:lin} with initial values $(\tilde E^{(0)}_n, \tilde H^{(0)}_n)$ such that
	\begin{equation}\label{eq:contraposition}
		\int_0^T \|\nu \cdot (\epsilon \tilde E_n, \mu \tilde H_n)\|_{L^2(\Gamma)}^2 \dd \tau
		> n \int_0^T \norm{\nu\times \tilde E_n}_{L^2(\Gamma)}^2 \dd \tau\,.
	\end{equation}
	We now divide $(\tilde E^{(0)}_n, \tilde H^{(0)}_n)$ by 
	\[
		\int_0^T \Big(\delta\norm{(\tilde E_n,\tilde H_n)}_{H^{-\theta}}^2 + c_\delta \norm{\nu \cdot (\epsilon \tilde E_n, \mu \tilde H_n)}_{H^{-1}(\Gamma)}^2\Big)\dd \tau\,.
	\]
	Using \eqref{eq:proof-time-independent-bnd} and \eqref{eq:contraposition}, we obtain classical solutions $(E_n, H_n)$ of \eqref{pde:lin} satisfying
	\begin{align}\label{eq:normal-traces-bounded}
		\int_0^T \norm{\nu \cdot (\epsilon E_n, \mu H_n)}_{L^2(\Gamma)}^2\dd \tau
		 & \leq \int_0^T \big(\delta\norm{(E_n,H_n)}_{H^{-\theta}}^2
		+  c\norm{\nu\times E_n}_{L^2(\Gamma)}^2                                                               \notag \\
		 & \hphantom{={}}+ c_\delta \norm{\nu \cdot (\epsilon E_n, \mu H_n)}_{H^{-1}(\Gamma)}^2\big)\dd \tau  \notag \\
		 & = 1 + c\int_0^T \norm{\nu\times E_n}_{L^2(\Gamma)}^2 \dd \tau                                       \notag \\
		 & \leq 1 + \frac cn \int_0^T \norm{\nu \cdot (\epsilon E_n, \mu H_n)}_{L^2(\Gamma)}^2 \dd \tau\,.
	\end{align}
	We thus derive the formulas
	\begin{align}\label{eq:proof-norm-trace-bdd}
		\int_0^T \left(\norm{\nu \cdot (\epsilon E_n, \mu H_n)}_{L^2(\Gamma)}^2\right) \dd \tau &\leq 2\,,\\
		\label{eq:proof-contradiction}
		\int_0^T \big(\delta\norm{(E_n,H_n)}_{H^{-\theta}}^2
		+ c_\delta \norm{\nu \cdot (\epsilon E_n, \mu H_n)}_{H^{-1}(\Gamma)}^2\big)\dd \tau &= 1\,,
	\end{align}
	for all sufficiently large $n \in \N$, and
	\begin{equation}\label{eq:proof-tang-trace-vanish}
		\int_0^T  \norm{\nu\times (E_n, H_n)}_{L^2(\Gamma)}^2  \dd \tau \to 0 \quad \text{ for $n\to\infty$, }
	\end{equation}
	where we also use \eqref{eq:contraposition}, the boundary condition in \eqref{pde:lin} and our construction.

	Next, Lemmas~\ref{lem:observability} and \ref{lem:SplitNormTang}, as well as formulas~\eqref{eq:proof-norm-trace-bdd} and \eqref{eq:proof-tang-trace-vanish} imply
	\[
		\norm{(E^{(0)}_n,H^{(0)}_n)}_{L^2} = \norm{(E^{(0)}_n,H^{(0)}_n)}_{H(\div_\epsilon) \times H(\div_\mu)}\lesssim 1 \quad \text{for }n\in\N\,.
	\]
	For that reason, $(E^{(0)}_n,H^{(0)}_n)$ converges weakly in $L^2(\Omega)$ to some limit $(E^{(0)},H^{(0)})$.
	(Here and below we pass to subsequences without mentioning it.)
	Note that the limit still satisfies $\div(\epsilon E^{(0)}) = \div(\mu H^{(0)}) = 0$.
	We denote the solution for $(E^{(0)}, H^{(0)})$ with $(E,H) \in C([0,T], L^2(\Omega))$.
	By continuity of the semigroup and the Rellich--Kondrachov theorem we obtain the convergence $(E_n,H_n)(t) \to (E,H)(t)$ in $H^{-\theta}$ for $t \in [0,T]$.
	Because of $(E_n,H_n),(E,H) \in \ker(\div_\epsilon) \times \ker(\div_\mu)$, we also have weak convergence in $H(\div_\epsilon) \times H(\div_\mu)$, and thus the normal trace converges weakly in $H^{-1/2}(\Gamma)$.
	Therefore, $\tr_n (E_n,H_n)$ tends to $\tr_n (E,H)$ strongly in $H^{-1}(\Gamma)$ for $n \to \infty$.

	From \eqref{eq:proof-contradiction} we conclude that $(E,H)$ also satisfies
	\begin{equation*}
		\int_0^T \big(\delta\norm{(E,H)}_{H^{-\theta}}^2
		+ c_\delta \norm{\nu \cdot (\epsilon E, \mu H)}_{H^{-1}(\Gamma)}^2\big)\dd \tau = 1\,.
	\end{equation*}
	In particular $(E,H) \neq 0$, which will lead to a contradiction.

	2.b) Lemma~3.2 of \cite{pokojovy2020} shows that
	$(E_n \times \nu , H_n \times \nu) \rightharpoonup (E \times \nu, H \times \nu)$
	weakly in
	$L^2([0,T] \times \Gamma)$ for $n \to \infty$ (see also Lemma 7.2.2.2 in \cite{lagnese2004}).
	Hence, $(E \times \nu, H \times \nu) = 0$ by \eqref{eq:proof-tang-trace-vanish}.

	We next regularize the fields $(E,H)$ by means of a Friedrich's mollifier as done in \cite{eller2007}.
	So we take $0 \leq \rho_{\tau_0} = \rho \in C^\infty_c(-\tau_0,\tau_0)$ with $\int_{-\tau_0}^{\tau_0} \rho \dd x =1$ and consider
	\begin{equation*}	\begin{aligned}
			E^*(t,x) & \coloneqq [\rho \ast E(\cdot,x)](t) = \int_0^T \rho(\tau)E(t-\tau,x) \dd \tau\,, \\
			H^*(t,x) & \coloneqq [\rho \ast H(\cdot,x)](t) = \int_0^T \rho(\tau)H(t-\tau,x) \dd \tau\,,
		\end{aligned} \qquad t \in [\tau_0, T-\tau_0]\,.
	\end{equation*}

	These are functions $(E^*, H^*)\in C^\infty\left([\tau_0, T-\tau_0],(L^2(\Omega))^6\right)$ which solve the partial differential equations in \eqref{pde:lin} and satisfy $(E^* \times \nu, H^* \times \nu) = 0$. %for $t \in [0,\infty)$.
	In particular, $(E^*, H^*)(t)$ belongs to $H(\curl) \cap X$.
	Applying $\partial_t$, we see that $(\partial_tE^*, \partial_tH^*)$ fulfill the same properties.
	Thus $(E^*, H^*)(t)$ and $(\partial_t E^*, \partial_t H^*)(t)$ belong to $H^1(\Omega)$ for $t \in [\tau_0, T-\tau_0]$ by Theorem~\ref{thm:div-curl-est}~a).

	Lemma~\ref{Lem:surfaceCurl} now implies that $\nu \cdot \curl E^* = \nu \cdot \curl H^* = 0$ on $\Gamma$.
	Notice that therefore $(\partial_t E^*, \partial_t H^*)$ are classical solutions of \eqref{pde:lin} with  vanishing normal traces $\nu \cdot \epsilon \partial_t E^* = \nu \cdot \mu \partial_t H^* = 0$.
	Lemmas~\ref{lem:SplitNormTang} and \ref{lem:observability} then show that $(\partial_t E^*, \partial_t H^*)(\tau_0) = 0$ and hence $(\partial_t E^*, \partial_t H^*)$ stays zero for $t \in [\tau_0, T-\tau_0]$ and so $(\curl E^*, \curl H^*)(t) = 0$ by \eqref{pde:lin}.
	Utilizing Theorem~\ref{thm:div-curl-est}~b) we see that $(E^*, H^*)$ vanishes on $[\tau_0, T-\tau_0]$.
	Letting $\tau_0 \to 0$, it follows that $(E^{(0)},H^{(0)}) = 0$, contradicting $(E,H) \neq 0$.
	Hence, \eqref{eq:proof-exp-stab:normal-trace} holds.

	3) For arbitrary $\tilde T = nT + t' \geq \hat T$ for some $n \in \N$ and $t'\in[0,T)$, using the translation invariance of the problem \eqref{pde:lin}, we calculate
	\begin{align} \label{eq:proof-exp-stab:normal-trace-multiples}
		\notag  \int_0^{\tilde T} \norm{\nu \cdot (\epsilon E, \mu H)}_{L^2(\Gamma)}^2 \dd \tau &= \sum_{k=1}^n\int_{(k-1)T}^{kT} \norm{\nu \cdot (\epsilon E, \mu H)}_{L^2(\Gamma)}^2 \dd \tau \\
		\notag &\quad+ \int_{nT}^{nT+t'} \norm{\nu \cdot (\epsilon E, \mu H)}_{L^2(\Gamma)}^2 \dd \tau \quad \\
		\notag &\lesssim \sum_{k=1}^n\int_{(k-1)T}^{kT} \norm{\nu\times E}_{L^2(\Gamma)}^2 \dd \tau + \int_{(n-1)T}^{nT+t'}\!\!\!\norm{\nu\times E}_{L^2(\Gamma)}^2 \dd \tau\\
		&\lesssim \int_0^{\tilde T} \norm{\nu\times E}_{L^2(\Gamma)}^2 \dd \tau\,.
	\end{align}
	For general $t\geq s\geq 0$ the result follows in the same way.\qedhere
\end{proof}

Lemmas~\ref{lem:SplitNormTang} and \ref{lem:observability}, the boundary condition in \eqref{pde:lin}, and Theorem~\ref{thm:normal-trace-vs-tang-trace} yield the observability estimate
\begin{equation}\label{eq:proof-time-independent-bnd-obs}
	\norm{(E,H)(0)}_{L^2}^2 \lesssim \int_0^{t}\norm{E \times \nu}_{L^2(\Gamma)}^2 \dd \tau
\end{equation}
for times $t \geq \hat T$, where $\hat T$ is given by Lemma~\ref{lem:observability}.

We further deduce the exponential stability of the semigroup $(S(t))_{t\geq0}$ on $X$. For scalar coefficients this fact was shown in \cite{eller2002-1}, \cite{eller2002-2}, and \cite{nicaise2003}. For the case of internal damping we refer to \cite{eller2019}, \cite{nicaise2005}, or \cite{phung2000}.
We note that the next result does not follow from Theorem~\ref{thm:mainThm}.

\begin{Corollary}
	Let \eqref{eq:cond-on-coefficients}, \eqref{eq:cond-on-initial-values}, \eqref{eq:non-trapping}, and \eqref{eq:complement-connected} hold. For $(E_0, H_0) \in X$ we have
	\[
		\norm{S(t)(E^{(0)}, H^{(0)})}_{L^2(\Omega)} \lesssim \e^{-\omega t}\norm{(E^{(0)}, H^{(0)})}_{L^2(\Omega)}
	\]
	for some $\omega>0$ and all $t \geq 0$.
\end{Corollary}

\begin{proof}
	Using Lemma~\ref{lem:InnerEnergyEst}, Lemma~\ref{lem:SplitNormTang}, Theorem~\ref{thm:normal-trace-vs-tang-trace}, as well as the boundary condition in \eqref{pde:lin} and \eqref{eq:energy-est-lin}, we obtain the estimate
	\begin{align*}
		\norm{(E,H)(t)}^2_{L^2} & + \int_0^{t} \int_\Omega (\epsilon E \cdot E + \mu H \cdot H )\dd x \dd \tau                                                    \\
		                         & \lesssim \norm{(E,H)(t)}^2_{L^2} + \norm{(E,H)(0)}^2_{L^2}                                                                      \\
		                         & \quad + \int_0^{t} \int_\Gamma\left( \abs{\nu \cdot(\epsilon E,\mu H)}^2 + \abs{\nu \times (E,H)}^2 \right) \dd \sigma \dd \tau \\
		                         & \lesssim \norm{(E,H)(t)}^2_{L^2} + \norm{(E,H)(0)}^2_{L^2} + \int_0^{t} \norm{E \times \nu}_{L^2(\Gamma)}^2 \dd \tau           \\
		                         & \lesssim \norm{(E,H)(0)}^2_{L^2}\,.
	\end{align*}
	for $t\geq \hat T$, where $\hat T$ is given by Lemma~\ref{lem:observability}.
	As in Remark~5.2 of \cite{lasiecka2019}, the exponential decay now follows in a standard way.
\end{proof}

Combining Theorem~\ref{thm:normal-trace-vs-tang-trace} with \eqref{eq:energy-est-lin} and Lemma~\ref{lem:SplitNormTang} we derive the following result, allowing us to bound the traces of the fields by their initial values globally in time.
(See Theorem~1.2 of \cite{cagnol2011} or Proposition~\ref{prop:gronwall} for results locally in time.)

\begin{Corollary}
	Let \eqref{eq:cond-on-coefficients}, \eqref{eq:cond-on-initial-values}, \eqref{eq:non-trapping}, and \eqref{eq:complement-connected} hold.
	Then the solution $(E,H) \in C([0,\infty), (L^2(\Omega))^6)$ of \eqref{pde:lin} satisfies
	\begin{equation*}
		\begin{split}
			\norm{(E(t), H(t))}^2_{L^2(\Omega)} + \int_s^t \norm{\tr(E,H)}^2_{L^2(\Gamma)} \dd \tau
			&\lesssim \norm{(E(s), H(s))}^2_{L^2(\Omega)}
		\end{split}
	\end{equation*}
	for $0\leq s \leq t$.
\end{Corollary}

\section{Proof of Proposition \ref{prop:regularity-boost}}
\label{sec:7}

We recall the crucial regularity result from Section~\ref{sec:3}.

\edef\tempThmCnt{\arabic{Theorem}}
\edef\tempSecCnt{\arabic{section}}
\setcounter{Theorem}{\regBoostThmCnt}
\setcounter{section}{\regBoostSecCnt}

\begin{Proposition}
	Under the conditions of Theorem \ref{thm:mainThm} with the exception of \eqref{eq:technicalCond} and the connectedness of $\R^3\setminus\Omega$ the following estimate holds
	\[
		z(t) \lesssim e(t)+z(t)^2
	\]
	for all $t \in [0,T_*)$.
\end{Proposition}

\setcounter{Theorem}{\tempThmCnt}
\setcounter{section}{\tempSecCnt}

Due to the nonlinear boundary conditions, one has to account for extra terms compared to the proof of Proposition~4.1 in \cite{pokojovy2020}, but we can still follow the reasoning given there.
Therefore, we will only outline the proof and focus on the differences to \cite{pokojovy2020}.
\begin{proof}
	Remember that $z$ is given by
	\[
		z(t) = \max_{0\leq j \leq 3} \left(\norm{\partial_t^jE(t)}_{H^{3-j}(\Omega)}^2 +\norm{\partial_t^j H(t)}_{H^{3-j}(\Omega)}^2 \right)
	\]
	and $e$ is defined as
	\[
		e(t) = \tfrac12 \max_{0\leq j \leq 3} \left(\norm{\widehat\epsilon_j^{1/2} \partial_t^jE(t)}_{L^2(\Omega)}^2 +\norm{\widehat\mu_j^{1/2} \partial_t^j H(t)}_{L^2(\Omega)}^2 \right)\,.
	\]
	The $L^2$-terms of $z$ can thus be trivially bounded by $e(t)$.
	The other (squared) norms will be estimated by $e(t)$ or $z(t)^2$ or by means of previous steps.

	Compared to \cite{pokojovy2020} we only have to modify the arguments that involve boundary conditions of differentiated problems.
	Here additional error terms appear as $\lambda$ also depends on $E$.
	They enter our reasoning only through the div-curl Theorem~\ref{thm:div-curl-est}~a) which controls the $H^1$-norm for fields $(u,v)$ by the estimate~\eqref{eq:appendiXA}.
	Here we have to control the $H^{1/2}$-norm of the boundary inhomogeneities on the right-hand side.
	We follow the steps of the proof given in Section~6 of \cite{pokojovy2020}.

	\textit{1) $H^1$-estimates for $\partial_t^k E$ and $\partial_t^k H$.}
	The (squared) $H^1$-norm of $\partial_t^k(E,H)$ for $k\in \{0,1,2\}$ is bounded via the div-curl estimate.
	Only for $k=2$ the boundary inhomogeneity $h_2$ is non-zero, where we have $\norm{h_2}_{H^{1/2}(\Gamma)} \lesssim z(t)$ by \eqref{eq:commutator-estimates}, as desired.

	\textit{2) Interior spatio-temporal estimates for $E$ and $H$.}
	One decomposes $(E,H) = (1-\chi) (E,H) + \chi(E,H)$ for a cutoff function $\chi$, being equal to $1$ near the boundary $\Gamma$. In the estimates for $(1-\chi) \partial_t^k (E,H)$ the boundary does not play a role so that the reasoning from \cite{pokojovy2020} has not to be modified.

	\textit{3) Preparation for the boundary collar estimates for $E$ and $H$.}
	To get $H^2$- or $H^3$- bounds on $\partial_t^k \chi (E,H)$, we employ tangential and normal derivatives $\partial_\tau$ and $\partial_\nu$ (extended to a neighborhood of $\Gamma$).
	We use the div-curl Theorem~\ref{thm:div-curl-est}~a) for $\partial_\tau \partial^k_t (\chi(E,H))$ with $k\in\{0,1\}$ and for $\partial_\tau^2 (\chi(E,H))$.
	Differentiating \eqref{eq:inhomogenous_maxwell_boundary} with $h_0=0$ or $h_1=0$ we obtain
	\begin{equation}\label{eq:QuasilinearBndCondition}
		\begin{aligned}
				&\partial_{\tau}\partial_{t}^{k}\chi H\times\nu+\left(\hat\lambda_k(\partial_{\tau}\partial_{t}^{k}\chi E\times\nu)\right)\times\nu
				= -\partial_{t}^{k}\chi H\times\partial_{\tau}\nu\\
				& \qquad-\left(\partial_{\tau}\hat\lambda_k(\partial_{t}^{k}\chi E\times\nu)\right)\times\nu
		-\left(\hat\lambda_k(\partial_{t}^{k}\chi E\times\partial_{\tau}\nu)\right)\times\nu\\
		&\qquad-\left(\hat\lambda_k(\partial_{t}^{k}\chi E\times\nu)\right)\times\partial_{\tau}\nu\,.
		\end{aligned}
	\end{equation}
	and 
	\begin{align*}
		&\partial_\tau^2(\chi H)\times \nu + \big(\lambda \partial_\tau^2(\chi E) \times \nu\big)\times \nu = - 2\partial_\tau(\chi H) \times \partial_\tau\nu - (\chi H) \times \partial^2_\tau\nu \\
		& \hphantom{{}=} - \big((\partial_\tau^2\lambda) (\chi E \times \nu)\big)\times \nu -\big(\lambda (\chi E) \times \partial_\tau^2\nu\big)\times \nu -\big(\lambda (\chi E) \times \nu\big)\times \partial_\tau^2\nu \\
		& \hphantom{{}=} 
		-2 \big((\partial_\tau\lambda) (\partial_\tau(\chi E) \times \nu)\big) \times \nu 
		- 2\big((\partial_\tau\lambda) (\chi E \times \partial_\tau\nu)\big)\times \nu\\ 
		& \hphantom{{}=} 
		- 2\big((\partial_\tau\lambda) (\chi E \times \nu)\big)\times \partial_\tau\nu 
		-2 \big(\lambda \partial_\tau(\chi E) \times \partial_\tau\nu\big)\times \nu \\
		& \hphantom{{}=} 
		- 2\big(\lambda (\partial_\tau\chi E) \times \nu\big)\times \partial_\tau\nu 
		- 2\big(\lambda \chi E \times \partial_\tau\nu\big)\times \partial_\tau\nu\,.
	\end{align*}
	with $\hat\lambda_0 = \lambda(\cdot, E)$ and $\hat\lambda_1 = \lambdad(\cdot, E)$.
	The remaining estimates involving $\partial_\nu$ are entirely based on the curl and div equations \eqref{eq:inhomogenous_maxwell} respectively \eqref{eq:divergence_inhomogeneities} and thus carry over from \cite{pokojovy2020}.

	\textit{4) $H^2$-estimates for $E$ and $H$.}
	When estimating $(E,H)$ in $H^2$, we need the boundary condition \eqref{eq:QuasilinearBndCondition} with $k=0$.
	Compared to \cite{pokojovy2020}, the only new term is 
	\[
		\tilde h_k \coloneqq \left(\partial_{\tau}\hat\lambda_k(\partial_{t}^{k}E\times\nu)\right)\times\nu
	\]
	with $k=0$, where we note that
	\begin{equation}\label{eq:differentiated-bnd-factor}
		\partial_j \hat\lambda_k(x,t) = (\partial_j \hat\lambda_k)(x,E(x,t)) + \sum_{i=1}^3 \partial_{\xi_i} \hat\lambda_k(x,E(x,t))\partial_jE_i(x,t)\,.
	\end{equation}
	We can thus bound this error term by
	\begin{equation}\label{eq:new-error-term}
	\begin{aligned}
		\bigl\|\tilde h_0\bigr\|_{H^{1/2}(\Gamma)} &\lesssim \norm{\partial_\tau \lambda(\cdot, E) E}_{H^1}\\
		&\lesssim \norm{\partial_x \lambda}_{W^{1,\infty}} \norm{E}_{H^1} + \norm{\partial_\xi \lambda}_{W^{1,\infty}} \norm{E}_{W^{1,\infty}} \norm{E}_{H^1}\\
		& \lesssim\norm{E}_{H^1} + z(t)
	\end{aligned}		
\end{equation}
as desired, using Sobolev's embedding and $z \leq 1$ by \eqref{eq:bound}.
(Observe that $\norm{E}_{H^1}$ and $\norm{\partial_t E}_{H^1}$ were already handled in step~\textit{1)}.)

	\textit{5) $H^2$-estimates for $\partial_t E$ and $\partial_t H$.}
	To bound $\partial_t(E,H)$ in $H^2$ we use \eqref{eq:QuasilinearBndCondition} with $k=1$. In \eqref{eq:new-error-term} one only has to replace $\norm{E}_{H^1}$ by $\norm{\partial_t E}_{H^1}$ in order to show
	\[
		\bigl\|\tilde h_1\bigr\|_{H^{1/2}(\Gamma)} \lesssim \norm{\partial_t E}_{H^1} +z(t)
	\]
	as before.

	\textit{6) $H^3$-estimates for $E$ and $H$.}
	We finally treat the $H^3$-norm of $(E,H)$.
	Compared to steps~\textit{4)} and \textit{5)}, in the boundary condition only the term $\hat h$ with $\partial_\tau^2 \lambda$ poses new difficulties.
	To tackle it, we differentiate \eqref{eq:differentiated-bnd-factor} once more in $x$ and employ Sobolev's embedding.
	It follows
	\begin{align*}
		\bigl\|\hat h\bigr\|_{H^{1/2}(\Gamma)} 
		&\lesssim \norm{\partial_\tau^2 \lambda(\cdot,E) E}_{H^1} 
		\lesssim \norm{\partial_x^2 \lambda}_{W^{1,\infty}} \norm{E}_{H^1} +\norm{\partial_{x\xi} \lambda}_{W^{1,\infty}} \norm{E}_{H^2}^2 \\
		&\hphantom{{}=} + \norm{\partial_{\xi}^2 \lambda}_{W^{1,\infty}} \norm{E}^2_{H^2} \norm{E}_{H^3} + \norm{\partial_{\xi} \lambda}_{W^{1,\infty}} \norm{E}_{H^3} \norm{E}_{H^2}\\
		&\lesssim\norm{E}_{H^1} + z(t)\,. \qedhere
	\end{align*}
\end{proof}

\appendix

\section{Div-curl estimate}

In this section we show the div-curl estimates in Theorem~\ref{thm:div-curl-est}.
To this aim, we have to study the Helmholtz decomposition in $H^{-\theta}$ to some extent.

\begin{Remark}
	Let $-\frac12 < s <\frac12$.
	Here the dual $(H^s(\Omega))^*$ equals $H^{-s}(\Omega)$ since $C_c^\infty(\Omega)$ is dense in $H^s(\Omega)$, see Theorem 11.1 in \cite{lions1972} or Chapter 4.3. in \cite{triebel1998}.
	As in \cite{fujiwara2007}, we define the normal trace $\tr_n \colon H^s(\div, \Omega) \to H^{s-\frac12}(\Gamma)$ distributionally by
	\[
		\langle \varphi, \tr_n u \rangle = \langle \div u, \Phi\rangle_{H^{s}(\Omega)\times H^{-s}(\Omega)} + \langle u, \nabla \Phi \rangle_{H^{s}(\Omega)\times H^{-s}(\Omega)}
	\]
	for $\varphi \in C^{\infty}(\Gamma)$ and $\Phi \in C_c^\infty(U)$ with $\tr \Phi = \varphi$ for a bounded domain $U \supset \overline\Omega$.
	This definition is independent of the choice of the continuation $\Phi$ of $\varphi$. Futhermore Theorem~2.3 in \cite{fujiwara2007} yields the continuity of $\tr_n$.

	Analogously we can define the tangential trace $\tr_t\colon H^s(\curl,\Omega) \to H^{s-\frac12}(\Gamma)$ by
	\[
		\langle \varphi, \tr_t u \rangle =  \langle u, \curl \Phi \rangle_{H^{s}(\Omega)\times H^{-s}(\Omega)} - \langle \curl u ,\Phi \rangle_{H^{s}(\Omega)\times H^{-s}(\Omega)}\,.
	\]
	We estimate
	\begin{align*}
		\langle \varphi, \tr_t u \rangle & \lesssim \norm{\Phi}_{H^{-s}(\Omega)}\norm{\curl u}_{H^{s}(\Omega)} + \norm{u}_{H^{s}(\Omega)}\norm{\curl \Phi}_{H^{-s}(\Omega)} \\
		                                 & \lesssim \norm{u}_{H^s(\curl)}\norm{\Phi}_{H^{1-s}(\Omega)}\,.
	\end{align*}
	Note that $\norm{\varphi}_{H^{\frac12-s}(\Gamma)}$ is the infimum of $\norm{\Phi}_{H^{1-s}(\Omega)}$, where $\Phi \in H^{1-s}(\Omega)$ with $\tr{\Phi}=\varphi$.
	Taking the infimum on both sides of the inequality above, we thus obtain
	\begin{align*}
		\langle \varphi, \tr_t u \rangle & \lesssim \norm{u}_{H^s(\curl)}\norm{\varphi}_{H^{\frac12-s}(\Gamma)}\,.
	\end{align*}
	This yields the continuity of $\tr_t\colon H^s(\curl) \to H^{s-1/2}(\Gamma)$ of the tangential trace.
	% because of
	% \[
	% 	\norm{\tr_t u}_{H^{s-\frac12}(\Gamma)} = \sup_{\varphi \in H^{\frac12-s}(\Gamma)} \frac{\langle \varphi, \tr_t u \rangle}{\norm{\varphi}_{H^{\frac12-s}(\Gamma)}} \lesssim \norm{u}_{H^s(\curl)}\,.
	% \]
\end{Remark}

Introducing
\begin{align*}
	H^s_{n0}(\div 0 , \Omega) = H^s_{n0}(\div 0) & \coloneqq \big\{f \in \big(H^s(\Omega)\big)^3 \mid \div f =0 ,\, \nu \cdot \tr f = 0 \big\}\,,                      \\
	G^s = G^s(\Omega)                            & \coloneqq \big\{u \in \big(H^s(\Omega)\big)^3 \mid \exists \varphi \in H^{s+1}(\Omega): \nabla \varphi = u\big\}\,,
\end{align*}
we can recall the following Helmholtz decomposition from Theorem~3.1 in \cite{fujiwara2007}.
\begin{Lemma} \label{lem:helmholtz}
	Let $\Omega \subseteq \R^3$ be bounded with $C^{2,1}$ boundary $\Gamma$. Further, suppose that $-\frac12<s<\frac 12$. We then have the following topological direct sum
	\[
		\big(H^s(\Omega)\big)^3 = H^s_{n0}(\div 0) \oplus G^s(\Omega)\,.
	\]
\end{Lemma}

\subsection*{The curl operator on $H^{1-\theta}$}

We state the estimate~(2.5) from \cite{eller2007}.
It is later utilized in an argument involving Peetre's lemma to show that the curl operator in a suitable setting has closed image.

\begin{Lemma}\label{lem:peetre-ineq}
	Let $\Omega \subseteq \R^3$ be bounded with $C^2$-boundary, $\theta \in [0,\frac12)$, and $v \in H^{1-\theta}(\Omega)$ with $\tr_n v = 0$.
	We then have
	\[
		\norm{v}_{H^{1-\theta}} \lesssim \norm{\div u}_{H^{-\theta}} + \norm{\curl v}_{H^{-\theta}} + \norm{v}_{H^{-\theta}}\,.
	\]
\end{Lemma}

The following result is well know for $\theta = 0$, see Remark~IX.1.4 of \cite{dautray1990}, where variants without connectivity assumptions are proven.

\begin{Lemma}\label{lem:curl}
	Let $\Omega \subseteq \R^3$ be bounded and \textit{simply connected} with a connected $C^{2,1}$ boundary $\Gamma$.
	For $\theta \in (0,\frac12)$ the curl as an operator
	\[
		\curl \colon (H^{1-\theta}_{n0}(\div 0), \norm{\cdot}_{H^{1-\theta}}) \to (H^{-\theta}(\div 0),\norm{\cdot}_{H^{-\theta}})
	\]
	is an isomorphism, and we have $\curl(H^{1-\theta}) = H^{-\theta}(\div 0)$.
\end{Lemma}

\begin{proof}
	1) Let $V \coloneqq (H^{1-\theta}_{n0}(\div 0), \norm{\cdot}_{H^{1-\theta}})$ denote the domain of the curl. Obviously $\curl V$ is a subspace of $\curl H^{1-\theta}$.
	Let $w\in H^{1-\theta}$.
	Lemma~\ref{lem:helmholtz} for $s=0$ yields functions $v \in H^0_{n0}(\div 0)$ and $\varphi \in H^{1}$ such that $w = v + \nabla \varphi$.
	It follows
	\[
		\Delta \varphi = \div \nabla \varphi = \div w \in H^{-\theta} \text{ with } \tr_n \nabla \varphi = \tr_n w \in H^{1/2 - \theta}(\Gamma)\,.
	\]
	Elliptic regularity implies that $\varphi$ belongs to $H^{2-\theta}$, so that $v = w - \nabla \varphi \in H^{1-\theta}$.
	As a result, $v$ is contained in $V$.
	Since $\curl w = \curl v$, the images $\curl H^{1-\theta} = \curl V$ agree.

	2) Next, we show that $\curl V$ is closed by means of Peetre's lemma (see Lemma~IX.1.2 in \cite{dautray1990}).
	To this end, we consider the operators $A_1 = \curl \colon V \to H^{-\theta}$ and $A_2 = \operatorname{Id} \colon V \to H^{-\theta}$.
	The Rellich--Kondrachov theorem shows that $A_2$ is compact.
	Lemma~\ref{lem:peetre-ineq} implies the inequality
	\[
		\norm{v}_{H^{1-\theta}} \lesssim \norm{A_1 v}_{H^{-\theta}} + \norm{A_2 v}_{H^{-\theta}}\quad \text{for } v\in V.
	\]
	Hence, $\curl V$ is closed in $H^{-\theta}$ by Peetre's lemma and, moreover, the kernel of curl on $V$ is finite dimensional.
	In the following step we show that it is indeed trivial.

	3) Let $w \in V$ with $\curl w = 0$ be an element of the kernel.
	In particular, $w$ belongs to $H(\curl 0)$.
	Thus, from Proposition~IX.1.2 of \cite{dautray1990} we obtain $w \in G^0 = \nabla H^1(\Omega)$, using simple connectedness.
	Since $H^{1-\theta}_{n0}(\div 0) \subseteq H^0_{n0}(\div 0)$, Lemma~\ref{lem:helmholtz} implies that $w=0$.
	Therefore, the kernel of $\curl\colon H^{1-\theta}_{n0}(\div 0) \to \curl(H^{1-\theta})$ is trivial.
	%Since the image is closed, the claim follows from the bounded inverse theorem.

	4) Finally, we establish that the image of the curl operator on $H^{1-\theta}(\Omega)$ is given by $\curl(H^{1-\theta}) = H^{-\theta}(\div 0)$.
	The result will then follow from the previous steps and the bounded inverse theorem.

	4.a) We begin by showing
	\begin{align}\label{eq:proof-curl:claim}
		\notag \curl(H^{1-\theta})^{\perp_\theta} & \coloneqq \{u \in H^{\theta}(\Omega) \mid \forall w \in \curl H^{1-\theta}(\Omega): \langle u , w\rangle_{H^{\theta}\times H^{-\theta}} = 0\} \\
		                                          & = H^{\theta}_{t0}(\curl 0) \coloneqq \{u \in H^\theta(\curl) \mid \curl u =0,\, \tr_t u = 0\}\,.
	\end{align}
	Let $u \in \curl(H^{1-\theta})^{\perp_\theta}$.
	For $v \in C_c^\infty (\Omega) \subseteq H^{1-\theta}(\Omega)$ we compute
	\[
		0 = \langle u , \curl v\rangle_{H^{\theta}\times H^{-\theta}} = \int_\Omega u \curl v \dd x = \langle  v, \curl u \rangle_{H_0^{1-\theta} \times H^{\theta- 1}}\,.
	\]
	Hence, $\curl u = 0$ and $u \in H^{\theta}(\curl)$.
	For $v \in H^1(\Omega)$ we then obtain
	\begin{align*}
		0 & = \langle u , \curl v\rangle_{H^{\theta}\times H^{-\theta}} = \langle u , \curl v\rangle_{H^{\theta}\times H^{-\theta}} - \langle \curl u, v \rangle_{H^{\theta}\times H^{-\theta}} \\
		  & = \langle\tr_t u , \tr v  \rangle_{H^{\theta-1/2}(\Gamma) \times H^{1/2-\theta})(\Gamma)},
	\end{align*}
	and the first inclusion $\curl(H^{1-\theta})^{\perp_\theta}\subseteq H^{\theta}_{t0}(\curl 0)$ follows.

	To show the reverse direction, take $u \in H^{\theta}_{t0}(\curl 0)$. Let $v \in H^{1-\theta}(\Omega)$ and $v_n\in C^\infty(\overline\Omega)\cap H^{1-\theta}(\Omega)$ with $v_n \to v$ in $H^{1-\theta}$.
	We have
	\[
		\langle u, \curl v_n\rangle_{H^{\theta}\times H^{-\theta}} = \langle \curl u ,v_n\rangle_{H^{\theta}\times H^{-\theta}} +\langle v_n, \tr_t u \rangle_{ H^{1/2-\theta}(\Gamma) \times H^{\theta-1/2}(\Gamma)} = 0\,.
	\]
	Since $\curl v_n \to \curl v$ in $H^{-\theta}$, we deduce $\langle u, \curl v\rangle_{H^{\theta}(\Omega)\times H^{-\theta}(\Omega)} = 0$. Therefore, $u$ belongs to $\curl(H^{1-\theta})^{\perp_\theta}$ and \eqref{eq:proof-curl:claim} is shown.

	4.b) We finish by proving
	\begin{equation}\label{eq:proof-curl:end}
		\curl(H^{1-\theta}) = H^{-\theta}(\div 0)\,.
	\end{equation}
	Let $u \in \curl(H^{1-\theta})$.
	Take $\varphi \in H^{1+\theta}_0$ and note that $\nabla \varphi \in H^\theta_{t0}(\curl 0)$, see for example the proof of Proposition IX.3 in \cite{dautray1990}.
	Step 4.a) thus yields
	\[
		0 = \langle  \nabla \varphi, u \rangle_{H^{\theta}(\Omega)\times H^{-\theta}(\Omega)} = \langle \varphi, \div u \rangle_{H^{1+\theta}(\Omega)\times H^{-1-\theta}(\Omega)}\,.
	\]
	Hence $\div u = 0$ and $u$ lies in $H^{-\theta}(\div 0)$.

	Next, let $u \in H^{-\theta}(\div 0)$ and take $v \in H^{\theta}_{t0}(\curl 0)$.
	Since $\Gamma$ is connected, Proposition IX.1.3 in \cite{dautray1990} yields a potential $\varphi \in H^{1+\theta}$ with $\nabla \varphi = v$ and constant trace $\tr \varphi = c$ for some $c \in \R$.
	By considering $\varphi -c$, we can assume that $\varphi \in H_0^{1+\theta}$.
	Take test functions $\varphi_k$ tending to $\varphi$ in $H^{1+\theta}$.
	We then calculate
	\begin{align*}
		\langle  v, u \rangle_{H^{\theta}(\Omega)\times H^{-\theta}(\Omega)}
		 & =\langle  \nabla \varphi, u \rangle_{H^{\theta}(\Omega)\times H^{-\theta}(\Omega)}                                                                                                           \\
		 & = \lim_{k \to \infty} \langle  \nabla \varphi_k, u \rangle_{H^{\theta}(\Omega)\times H^{-\theta}(\Omega)} + \langle \varphi_k, \div u \rangle_{H^{\theta}(\Omega)\times H^{-\theta}(\Omega)} \\
		 & = \lim_{k \to \infty} \langle \tr \varphi_k, \tr_n u \rangle = 0\,.
	\end{align*}
	We conclude that
	\begin{align*}
		u \in & \prescript{\perp_\theta}{}{(H^{\theta}_{t0}(\curl 0))} \coloneqq \{u \in H^{-\theta} \mid \forall v \in H^{\theta}_{t0}(\curl 0) : \langle v,u \rangle_{H^{\theta}\times H^{-\theta}} =0 \} \\
		      & = \prescript{\perp_\theta}{}{(\curl(H^{1-\theta})^{\perp_\theta})} = \overline{\curl(H^{1-\theta})} = \curl(H^{1-\theta})\,,
	\end{align*}
	using \eqref{eq:proof-curl:claim} and step 2). Therefore, \eqref{eq:proof-curl:end} is true.
\end{proof}

We next observe that Lemma~2.4 in \cite{cessenat1996} extends from $\theta =0$ to $\theta \in (-\frac12, 0]$.

	\begin{Lemma}\label{lem:cessenat2.4}
		Let $u\in H^{-\theta}(\curl)$ and $\theta\in[0,\frac12)$.
		Then
		\[
			\curl_\Gamma(\tr_\tau u) = -\div_\Gamma(\tr_t u) \in H^{-\theta - \frac12}\,,
		\]
		where we set
		\[
			\curl_\Gamma (\tr_\tau u) \coloneqq \tr_n(\curl u)
		\]
		and, for all $\varphi \in C^\infty(\Gamma)$ and $\Phi \in C^\infty(U)$ with $\tr \Phi = \varphi$ and open $U \supseteq \overline \Omega$,
		\[
			\langle\div_\Gamma(\tr_t u), \varphi\rangle \coloneqq - \langle\tr_t u , \tr_\tau(\nabla \Phi)\rangle\,.
		\]
	\end{Lemma}

	\begin{proof}
		First note that $\curl u \in H^{-\theta}(\div 0)$ and thus $\tr_n \curl u \in H^{-\theta-\frac12}(\Gamma)$.
		Since also $\tr_t u \in H^{-\theta- \frac12}(\Gamma)$, the definitions yield
		\begin{align*}
			\langle\curl_\Gamma(\tr_\tau u), \Phi\rangle & = \langle \tr_n(\curl u), \Phi\rangle =\langle \Phi, \div \curl u \rangle_\Omega + \langle \nabla \Phi, \curl u \rangle_\Omega                                  \\
			                                             & = - \langle\tr(\nabla\Phi),\tr_t u\rangle = - \langle \tr_\tau(\nabla \Phi), \tr_t u \rangle_{H^{\theta + \frac12}(\Gamma)\times H^{-\theta - \frac12}(\Gamma)}
			\\
			                                             & = - \langle \varphi, \div_\Gamma(\tr_t u) \rangle\,,
		\end{align*}
		where we partly omit the underlying spaces.
	\end{proof}

	\subsection*{The div-curl estimate}

	We show a corrected and improved version of Lem\-ma~5.1 of \cite{pokojovy2020}.
	We note that part~a) was shown in Lemma~4.5.5 of \cite{costabel2010} for scalar $\gamma$ and $\theta = 0$ in a similar way.
	In \cite{pokojovy2020} the $L^2$-norms of $u$ and $v$ on the right-hand side of \eqref{eq:appendiXA} were omitted erroneously.
	The estimate~\eqref{eq:appendiXA} for $\theta = 0$ below suffices for \cite{pokojovy2020}.
	Actually this version does not require additional geometric properties.
	In part~b) we can remove the $L^2$-norms on the right-hand side, as needed to obtain global-in-time properties in Lemma~\ref{lem:RandtermH-1k0}  and Theorem~\ref{thm:normal-trace-vs-tang-trace}, assuming that $\R^3\setminus\Omega$ is connected.
	This condition is needed, since without it the operator $A$ in the proof would have a non-trivial kernel of the form $(\nabla \varphi, \nabla \psi)$, where $\div(\alpha \nabla \varphi) = \div(\beta \nabla\psi) = 0$ and $\varphi$ and $\psi$ are constant on components of $\Gamma = \partial \Omega$.

\begin{Theorem}\label{thm:div-curl-est}
	Let $\Omega \subseteq \R^3$ be bounded with boundary $\Gamma\in C^2$ and $\theta \in [0, \frac12)$.
	Assume that $\alpha, \beta \in W^{1, \infty}(\Omega, \Rsym)$, $\gamma \in W^{1,\infty}(\Gamma, \Rsym)$ satisfy $\alpha, \beta, \gamma \geq \eta >0$ and $\gamma \nu^\bot \subseteq \nu^\bot$.
	Let $(u,v) \in H^{-\theta}(\curl)$ fulfill $\div(\alpha u) \in H^{-\theta}$, $\div(\beta v) \in H^{-\theta}$, and $v \times \nu + \gamma(u \times \nu) \times \nu \eqqcolon h \in H^{1/2 - \theta}(\Gamma)^3$.
	\begin{enumerate}[label=\alph*), wide, labelwidth=0pt, labelindent=0pt]
		\item Then $u$ and $v$ belong to $H^{1-\theta}(\Omega)^3$ and
		      \begin{equation}\label{eq:appendiXA}
			      \begin{aligned}
				      \norm{u}_{H^{1-\theta}} + \norm{v}_{H^{1-\theta}} & \lesssim
				      \norm{u}_{H^{-\theta}(\curl)} + \norm{v}_{H^{-\theta}(\curl)} + \norm{\div(\alpha u)}_{H^{-\theta}}                                                               \\
				                                                        & \quad+ \norm{\div(\beta v)}_{H^{-\theta}} + \norm{h}_{H^{1/2 - \theta}(\Gamma)} \eqqcolon N(u,v,h)\,.
			      \end{aligned}
		      \end{equation}
		\item Let $\theta = 0$. Assume in addition that $\R^3\setminus \Omega$ is connected. We then obtain
		      \begin{equation}\label{eq:appendiXB}
			      \begin{aligned}
				      \norm{u}_{H^{1}} + \norm{v}_{H^{1}} & \lesssim \norm{\curl u}_{L^2} + \norm{\curl v}_{L^2} + \norm{\div(\alpha u)}_{L^2} \\
				                                          & \quad + \norm{\div(\beta v)}_{L^2} + \norm{h}_{H^{1/2}(\Gamma)}\,.
			      \end{aligned}
		      \end{equation}
	\end{enumerate}
\end{Theorem}

\begin{proof}
	a) As in Lemma~4.5.5 of \cite{costabel2010} or Proposition~6.1 of \cite{lasiecka2019}, we use a finite partition of unity $\chi_j$ for $\overline\Omega$ such that $\supp \chi_j \subseteq \overline\Omega_j \subseteq \overline\Omega$ for  starshaped, bounded open sets $\Omega_j$.
	Following Lemma~5.1 of \cite{pokojovy2020} we will prove that $\chi_j u$ and $\chi_j v$ belong to $H^{1-\theta}(\Omega_j)$ and
	\[
		\norm{\chi_j u}_{H^{1-\theta}} + \norm{\chi_j v}_{H^{1-\theta}} \leq N(\chi_j u,\chi_j v,\chi_j h) \lesssim N( u,v, h)
	\]
	Summing these pieces, we obtain \eqref{eq:appendiXA}.

	So let $\Omega$ be starshaped.
	Since $\curl u \in H^{-\theta}(\div 0)$, by Lemma~\ref{lem:curl} there is a field $w \in H^{1-\theta}_{n0}(\div0)$ such that $\curl u = \curl w$ and $\norm{w}_{H^{1-\theta}} \lesssim \norm{\curl u}_{H^{-\theta}}$.
	So $u-w$ belongs to $H^{-\theta}(\curl 0)$.
	We show that $H^{-\theta}(\curl 0) = G^{-\theta}$.

	Indeed, the inclusion ``$\supseteq$'' is clear.
	Conversely, take $f \in H^{-\theta}(\curl0)$.
	Lemma~\ref{lem:helmholtz} provides maps $v_0 \in H^{-\theta}_{n0}(\div0)$ and $\varphi \in G^{-\theta}$ such that $f=v_0 + \nabla \varphi$.
	Hence, $v_0$ is contained in $H^{-\theta}_{n0}(\div0) \cap H^{-\theta}(\curl 0)$.
	For an open, simply connected set $U \subseteq \overline U \subseteq \Omega$, take a cut-off function $\chi \in C_c^\infty(\Omega)$ with $\chi=1$ on $U$ and let $\tilde v$ be the zero-extension of $\chi v_0$ to $\R^3$.
	Computing
	\begin{align*}
		\norm{\div \tilde v}_{H^{-\theta}(\R^3)}
		 & = \sup_{\varphi \in H^\theta(\R^3)} \abs{\frac{\langle \chi\nabla\varphi,v_0 \rangle}{\norm{\varphi}_{H^{\theta}}}}
		= \sup_{\varphi \in H^\theta(\R^3)} \frac{\abs{\langle \varphi\nabla\chi,v_0 \rangle}}{\norm{\varphi}_{H^{\theta}}}\\
		&\lesssim_\chi \sup_{\Psi \in (H^\theta)^3} \frac{\abs{\langle \Psi,v_0 \rangle}}{\norm{\Psi}_{H^{\theta}}}
		= \norm{v_0}_{H^{-\theta}}\,,                                                                                        \\
		\norm{\curl \tilde v}_{H^{-\theta}(\R^3)}
		 & = \sup_{\Phi \in H^\theta(\R^3)^3} \abs{\frac{\langle \chi\curl\Phi,v_0 \rangle}{\norm{\Phi}_{H^{\theta}}}}
		= \sup_{\Phi \in H^\theta(\R^3)^3} \frac{\abs{\langle \Phi \times \nabla \chi,v_0 \rangle}}{\norm{\Phi}_{H^{\theta}}}\\
		&\lesssim_\chi \sup_{\Psi \in (H^\theta)^3} \frac{\abs{\langle \Psi,v_0 \rangle}}{\norm{\Psi}_{H^{\theta}}}= \norm{v_0}_{H^{-\theta}}\,,
	\end{align*}
	we see that $\tilde v, \div \tilde v$ and $\curl \tilde v$ belong to $H^{-\theta}(\R^3)$.
	Hence, $\hat v \coloneqq (1-\Delta)^{-\theta/2}\tilde v$ is an element of $H(\div,\R^3) \cap H(\curl, \R^3) = H^1(\R^3)$.
	As a result, $\tilde v$ is contained in $H^{1-\theta}(\R^3)$ and $v_0$ in $H^{1-\theta}_{\mathrm{loc}}(\Omega)$, which means that $v_0 \in H(\curl 0, U)$.
	Proposition~IX.1.2 in \cite{dautray1990} now yields a map $\Psi_U \in H^1(U)$ with $\nabla \Psi_U = v_0$ on $U$.
	Because $\div v_0 = 0$ the map $\Psi_U$ is analytic, and by analytic continuation we find a smooth function $\Psi$ on $\Omega$ with $\nabla \Psi = v_0 \in H^{-\theta}$.

	Let $x_0$ be the star center of $\Omega$.
	We can write $\Psi = \Psi(x_0) + \int_0^1 T_\tau \nabla \Psi \cdot m \dd \tau$, where $T_\tau g(x) = g(x_0 - \tau (x-x_0))$ and $m(x) = x-x_0$.
	This identity implies that $\Psi$ belongs to $H^{-\theta}(\Omega)$.
	Since also $\Delta \Psi = 0$ on $\Omega$ and $\partial_\nu \Psi = \tr_n v_0 = 0$ on $\Gamma$, we deduce $\Psi=0$ and so $v_0 = 0$, as needed.

	Hence, there exists a function $\varphi \in H^{1-\theta}$ such that $\nabla \varphi = u -w$ and therefore
	\begin{equation}\label{eq:proof-div-curl}
		\div(\alpha \nabla \varphi) = \div(\alpha u) - \div(\alpha w) \in H^{-\theta}\,.
	\end{equation}
	Here we can assume that $\int_\Omega \varphi \dd x = 0$.
	We thus obtain 
	\[
		\norm{\varphi}_{H^{-\theta}} \lesssim \norm{\nabla \varphi}_{H^{-\theta}} \lesssim \norm{u}_{H^{-\theta}} + \norm{w}_{H^{-\theta}} \lesssim \norm{u}_{H^{-\theta}(\curl)}	
	\]
	by the Poincaré inequality in $H^{-\theta}(\Omega)$ which can be shown as in $L^2(\Omega)$. %follow proof in Evans PDE
	One also infers
	\[
		\norm{\tr \varphi}_{H^{-\theta}(\Gamma)} \lesssim \norm{\varphi}_{H^{1-\theta}} \lesssim \norm{u}_{H^{-\theta}(\curl)}\,.
	\]
	Let $f\times \nu = \begin{psmallmatrix}
			0 & \nu_3 & - \nu_2\\
			-\nu_3 & 0 & \nu_1\\
			\nu_2 & -\nu_1 &0
		\end{psmallmatrix} f \eqqcolon Bf$ and $\tilde \gamma = B^\intercal\gamma B$.
	As in Lemma~5.1 of \cite{pokojovy2020} we now see that 
	\begin{align*}
		\div_\Gamma(\tilde\gamma \tr_\tau(\nabla\varphi)) &= \div_\Gamma(\gamma(w \times \nu)\times \nu) - \div_\Gamma(\gamma(u \times \nu) \times \nu)\\
		&= \div_\Gamma(\gamma(w \times \nu)\times \nu) - \div_\Gamma(h) + \div_\Gamma(\tr_t v)\\
		&= \div_\Gamma(\gamma(w \times \nu) \times \nu) - \div_\Gamma(h) - \tr_n(\curl v)
	\end{align*}
	belongs to $H^{-\theta-\frac12}(\Gamma)$, using also Lemma~\ref{lem:cessenat2.4} and that $\curl v \in H^{-\theta}(\div0)$.
	Elliptic regularity then yields $\tr \varphi \in H^{\frac32-\theta}(\Gamma)$, and then $\varphi \in H^{2-\theta}$ because of \eqref{eq:proof-div-curl}.
	Collecting the above relations we can thus bound $u$ in $H^{1-\theta}$ as asserted.
	One can treat $v$ similarly also using $u\in H^{1-\theta}$ and the boundary condition.
	(Here one can let $\tilde \gamma = I$).

	b) Now assume that $\R^3 \setminus \Omega$ is connected and that $\theta = 0$.
	We show that the $L^2$-norms appearing on the right-hand side of \eqref{eq:appendiXA} can be estimated by the div-curl-terms.
	Assume that
	\[
		\norm{(u,v)}_{H^1} \not\lesssim \norm{\curl u}_{L^2} + \norm{\curl v}_{L^2} +\norm{\div (\alpha u)}_{L^2}+\norm{\div (\beta v)}_{L^2} + \norm{h}_{H^{1/2}(\Gamma)}\,.
	\]
	Hence, there exists a sequence $(u_k,v_k) \in H^1$ with $\norm{(u_k,v_k)}_{H^1} = 1$ such that the right-hand side tends to zero.
	By means of the Banach--Alaoglu and Rellich--Kondrachov theorems, we can choose a subsequence, again denoted by $(u_k,v_k) $, which converges to a limit $(\overline u,\overline v)$ in $L^2$.

	We employ the operator $A = \begin{psmallmatrix}0 &\curl \\ -\curl & 0\end{psmallmatrix}$ on $X = \left(\ker(\div_\alpha) \cap \ker(\div_\beta), \norm{\cdot}_{L^2}\right)$
	with
	\[
		D(A) \coloneqq \{(u,v) \in X \mid \curl u , \curl v \in L^2, v \times \nu + \gamma (u \times \nu) \times \nu = 0\}	\,.
	\]
	Note that $(\overline u,\overline v)$ belongs to the kernel of $A$.

	We will prove below that $A$ is injective.
	Therefore, the sequence $(u_k,v_k)$ converges to $(\overline u,\overline v) = 0$ in $L^2$. The estimate \eqref{eq:appendiXA} with $\theta = 0$ then shows that $\norm{(u_k,v_k)}_{H^1} \to 0$, contradicting $\norm{(u_k,v_k)}_{H^1} = 1$ and the assertion follows.

	So take $(u,v) \in D(A)$ with
	$0 = A\begin{psmallmatrix}
			u\\v
		\end{psmallmatrix}$
	and consider
	\begin{align*}
		0 & = \left(A\begin{psmallmatrix}
				             u\\v
			             \end{psmallmatrix} \,\Big|\, \begin{psmallmatrix}
				                                          u\\v
			                                          \end{psmallmatrix}\right) = \int_\Omega (\curl v \cdot u - \curl u \cdot v) \dd x\,.
		\intertext{Integrating by parts and using the boundary condition, we compute}
		0 & = \int_\Gamma (u\times \nu) \cdot v \dd \sigma = -\int_\Gamma u \cdot (v \times \nu) \dd \sigma
		= \int_\Gamma (u\times \nu)\cdot \gamma  (u\times \nu) \dd \sigma\,.
	\end{align*}
	The positive definiteness of $\gamma$ implies that $u \times \nu = 0$ and therefore $v \times \nu = 0$ as well.
	Thus, $(u,v) \in \ker A$ satisfies
	\begin{align*}
		\curl u = \curl v = 0\,,\quad &\text {on } \Omega,\\
		u \times \nu = v \times \nu = 0\,,\quad &\text {on } \Gamma.
	\end{align*}
	Theorem~2.8 of \cite{cessenat1996} yields functions $\varphi, \psi \in H^1(\Omega)$ with $\nabla \varphi = u$ and $\nabla \psi = v$, which are constant on each component of $\partial \Omega$.
	Because of $(u,v)\in X$ we also have $\div (\alpha \nabla \varphi) = 0$ and $\div (\beta \nabla \psi) = 0$.
	Since $\R^3\setminus\Omega$ is connected, we see that $\varphi$ and $\psi$ are constant on each component of $\Omega$, so that $u=0=v$. Hence, $A$ is injective as desired. \qedhere
\end{proof}

\section*{Acknowledgments}
Roland Schnaubelt thanks Irena Lasiecka and Michael Pokojovy for fruitful discussions.

Funded by the Deutsche Forschungsgemeinschaft (DFG, German Research Foundation) – Project-ID 258734477 – SFB 1173.

\bigskip

\textbf{Declarations of interest: none.}

\printbibliography
\end{document}